\documentclass[12pt,twoside]{article}
\usepackage{amsmath,amsbsy,amsfonts}
\usepackage{hyperref}
\usepackage{color}

\newtheorem{lem}{Lemma}[section]
\newtheorem{prop}{Proposition}[section]
\newtheorem{thm}{Theorem}[section]

\newtheorem{cor}{Corollary}[section]
\newtheorem{rmk}{Remark}[section]

\newcommand{\cqd}{\hspace{10pt}\fbox{}}

\let\Section=\section

\def\section{\setcounter{equation}{0}\Section}

\def\nd{\noindent}

\def\proof{{\rm \bf Proof}}

\newcommand{\w}{W_0^{1,\Phi}(\Omega)}

\def\r{\bf{R}}

\begin{document}

\title{Existence and regularity of positive solutions of quasilinear elliptic  problems with  singular semilinear term}

\author{Jos\'e V. A. Goncalves~~~Marcos L. M. Carvalho 
\\
Carlos Alberto Santos\footnote{ Supported by CAPES/Brazil Proc.  $N^o$ $2788/2015-02$,}}

\date{}

\pretolerance10000

\maketitle

\begin{abstract}
\small{
\noindent {\small This paper deals with existence and regularity of positive solutions  of  singular elliptic problems on a smooth bounded domain with Dirichlet boundary conditions involving the $\Phi$-Laplacian operator.  The proof of existence is based on a variant  of the generalized Galerkin method that we developed inspired on ideas by Browder \cite{Browder} and a comparison principle. By using a kind of Moser iteration scheme we show $L^{\infty}(\Omega)$-regularity for positive solutions.}}
\end{abstract}

\begin{center}
{\bf Contents}
\end{center}

{ 1. Introduction}

{ 2. Main Results}

{ 3. A family of Auxiliary Problems}

{ 4. Applied Generalized Galerkin Method}

{ 5. Comparison of Solutions and  Estimates}

{ 6. Proof of  the Main Results }

   \hspace{0.4 cm} {6.1. Pure Singular Problem - Existence of Solutions}
    
   \hspace{0.4 cm} {6.2. Convex Singular Problem - Regularity of Solutions}

{ 7. Appendix  - On  Orlicz-Sobolev spaces}

\section{Introduction}

This paper concerns existence and regularity of solutions  to the  singular elliptic problem 
\begin{equation}\label{prob}
\displaystyle -\mbox{div}(\phi(|\nabla u|)\nabla u) = \frac{a(x)}{u^{\alpha}}~\mbox{in}~\Omega,~~
u>0~\mbox{in}~\Omega,~u=0~\mbox{on}~\partial \Omega,
\end{equation}
where $\Omega\subset\mathrm{R}^N$, with $N\geq 2$, is a bounded domain with smooth boundary $\partial \Omega$, $a$ is a non-negative function, $0 < \alpha< \infty$ and   $\phi:(0,\infty) \rightarrow (0,\infty)$ is  of class  $C^1$ and satisfies
\begin{itemize}
	\item[($\phi_1$)]  \     \ $\mbox{(i)} \ \ t\phi(t)\to 0 \ \mbox{\it as} \  t\to 0, ~~ \mbox{(ii)} \  t\phi(t)\to\infty \ \mbox{\it as} \  t\to\infty$,
	\item[($\phi_2$)]  \ $t\phi(t) \ \mbox{\it is strictly increasing in}~ (0, \infty)$,
	\item[($\phi_3$)] {\it there exist $\ell,m\in(1,N)$ such that }
	$$\ell-1\leq \frac {(t\phi(t))^\prime}{\phi(t)}\leq m-1,~t>0.$$
\end{itemize}
\nd We extend   $s \mapsto s \phi(s)$ to $\mathrm{R}$ as an odd function. It follows that the function 
$$
\Phi(t)=\int_0^ts\phi(s)ds,~t\in\mathrm{R}
$$
\nd is even and it is actually an $N$-function. Due to the nature of the operator
$$
\Delta_\Phi u :=\mbox{div}(\phi(|\nabla u|)\nabla u)
$$
we shall work in the framework of Orlicz and Orlicz-Sobolev spaces namely $L_{\Phi}(\Omega), L_{\widetilde{\Phi}}(\Omega)$ and $\w$. 

\nd We recall some basic notation on these spaces along with bibliographycal references in the Apendix.

\nd In the last years many research papers have been devoted to the study of singular problems like (\ref{prob}).   In \cite{karlin}, Karlin \& Nirenberg studied the singular integral equation
$$
u(x) = \int_0^1 G(x,y) \frac{1}{u(y)^{\alpha}} dy,~0 \leq x \leq 1,
$$
where $\alpha >0$ and $G(x,y)$ is a suitable potential. In \cite{crandall-rab-tartar}, Crandall Rabinowitz   \&  Tartar, addressed a class of singular problems which included as a special case, the  model problem
\begin{equation}\label{model prob}
\displaystyle-\Delta  u=  \frac{a(x)}{u^{\alpha}}~\mbox{in}~\Omega,~~ 
  u>0~\mbox{in}~\Omega,~u=0~\mbox{on}~\partial \Omega,
  \end{equation}
\nd where $\alpha > 0$ and $a: \Omega  \to [0, \infty)$  is a suitable $L^{1}$-function.  A broad literature on problems like \eqref{model prob} is available to date. We  would like to mention \cite{LS, CD-1, ZCheng-1}  and their references. We would like to refer the reader to  the very recent papers by Orsina \& Petitta \cite{OP-1},  Canino,   Sciunzi  \&   Trombetta  \cite{canino-0}   for the problem
$$
 \displaystyle -\Delta u= \frac{\mu}{u^{\alpha}}~  \mbox{in}\ \Omega,~ u>0\ \mbox{in}\ \Omega,\ u= 0~\mbox{on}~\partial\Omega.
$$
\nd In \cite{OP-1}   $\mu$ is a nonnegative bounded Radon measure while in \cite{canino-0} $\mu$ is an $L^{1}$ function. Other kinds of operators have been addressed and we mention  Chu-Wenjie \cite{chu} and De Cave \cite{cave}  for problems involving the p-Laplacian like
$$
-{\rm div}(|\nabla u|^{p-2}\nabla u)= \frac{a(x)}{u^{\alpha}}~  \mbox{in}\ \Omega,~ u>0\ \mbox{in}\ \Omega,\ u= 0~\mbox{on}~\partial\Omega;$$
\nd Qihu Zhang  \cite{qihu} and  Liu, Zhang \& Zhao \cite{liu-0} for $p(x)$-Laplacian operator,
$$
-{\rm div}(|\nabla u|^{p(x)-2}\nabla u)= \frac{a(x)}{u^{\alpha}}\ \  \mbox{in}\ \Omega,~~ u>0\ \mbox{in}\ \Omega,\ u= 0~\mbox{on}~\partial\Omega;
$$
\nd Boccardo \& Orsina \cite{bocardo}  and Bocardo \& Casado-D\'iaz \cite{casado}  for the problem   
 $$
 -{\rm div}(M(x) \nabla u) = \frac{a(x)}{u^{\alpha}}\ \  \mbox{in}\ \Omega,~~ u>0\ \mbox{in}\ \Omega,\ u= 0~\mbox{on}~\partial\Omega,
$$
where $M$ is a suitable matrix, 
\nd Lazer \& McKeena \cite{lazerMckeena-2}; Goncalves \& Santos \cite{GS-1}, Hu \& Wang \cite{HuWang-1} for  problems involving the Monge-Amp\'ere operator, e. g., 
$$
\mbox{\rm det} (D^2 u) = \frac{a(x)}{(-u)^{\gamma}}~~ \mbox{in}~~ \Omega,~u < 0~\mbox{in}~\Omega,~ u = 0~\mbox{on}~\partial\Omega,
$$
\nd where $a \in C^{\infty}(\overline{\Omega})$, $a > 0$ and $\gamma > 1$.

\nd To the best of our knowledge singular problems like \eqref{prob} in the presence of the operator $\Delta_{\Phi}$ were never studied  and the main results of this paper (see Section \ref{MAIN res}) namely Theorems \ref{Teor-prin}, \ref{Teor-prin1} as well as Corollary \ref{corol-1}  are new.

\nd Other problems which are special cases of \eqref{prob} are
\begin{eqnarray}
	-\Delta_p u-\Delta_q u={a(x)}{u^{-\alpha}}~\mbox{in}~\Omega,~~
	u>0~\mbox{in}~\Omega,~u=0~\mbox{on}~\partial \Omega,
\end{eqnarray}
\nd where $\phi(t)=t^{p-2}+t^{q-2}$ with $1<p<q<N$,
\begin{eqnarray}
-\sum_{i=1}^N\Delta_{p_i} u= {a(x)}{u^{-\alpha}}~\mbox{in}~\Omega,~~ 
u>0~\mbox{in}~\Omega,~u=0~\mbox{on}~\partial \Omega.
\end{eqnarray}
\nd where $\phi(t) = \sum_{j=1}^{N} t^{p_{j}-2}$, $1 < p_{1} < p_{2} < \ldots < p_{N} < \infty$ and $\sum_{j=1}^{N} \dfrac{1}{p_{j}} > 1$, 
\begin{eqnarray}\label{ex-geral}
\displaystyle-\mbox{div}(a(|\nabla u|^p)|u|^{p-2}\nabla u) = { a(x)}{u^{-\alpha}}~\mbox{in}~\Omega,~~ 
u>0~\mbox{in}~\Omega,~u=0~\mbox{on}~\partial \Omega.
\end{eqnarray}
\nd where  $\phi(t)=a(t^p)t^{p-2}$, $2\leq p<N$ and $a: (0,\infty) \to (0, \infty)$ is a suitable  $C^1(\mathbb{R}^+)$-function.

\nd  We also refer the reader to the  paper  \cite{repovs}, where the operator $\Delta_\Phi$ is employed. The operator $\Delta_\Phi$  appears in  applied mathematics, for instance in Plasticity, see e.g. Fukagai  and Narukawa \mbox{\rm \cite{Fukagai}} and references therein. We refer the reader to \cite{radu-3} for problems involving  general operators.

\section{Main Results}\label{MAIN res}

\nd   In this work, for each $x\in\Omega$, we set $
\displaystyle d(x)=\inf_{y\in\partial \Omega}  |x-y|$. Our first result is.
\begin{thm}\label{Teor-prin}
{\it	Assume that $(\phi_1)-(\phi_3)$ and $a \in  L^1(\Omega)$ hold. 
\nd  Then there is $u$ such that $
u^{(\alpha-1+\ell)/\ell} \in W_0^{1,\ell}(\Omega)$, $u \geq Cd~\mbox{a.e. in}~\Omega$, for some $C>0$, and: 
\vskip.2cm
\nd  $\rm \bf {(i)}$   $u \in \w$, and 
\begin{eqnarray}\label{final2}
\int_\Omega\phi(|\nabla u|)\nabla u\nabla \varphi dx=\int_\Omega\frac{a(x)}{u^\alpha}\varphi dx,~ \varphi  \in \w ,
\end{eqnarray}
}
\nd provided additionally that either
$
\displaystyle {a}{d^{-\alpha}}\in L_{\widetilde\Phi}(\Omega)
$
or $0<\alpha \leq 1$ and $a \in  L^{\ell^*/(\ell^*+ \alpha - 1)}(\Omega)$,
\vskip.2cm
\nd $\rm \bf {(ii)}$    $u \in W^{1,\Phi}_{loc}(\Omega)$, 
\nd and
\begin{equation}\label{final1}
\int_\Omega\phi(|\nabla u|)\nabla u\nabla \varphi dx=\int_\Omega\frac{a(x)}{u^\alpha}\varphi dx,~ \varphi \in C_0^{\infty}(\Omega) 
\end{equation}
\nd provided in addition that  $\alpha \geq 1$.
\end{thm}
\vskip.2cm

\nd Next we will present some regularity  results:

\begin{cor}\label{corol-1}
 Under the conditions of the above Theorem, we have that:
\begin{enumerate}
\item [$(i)$] $u \in C(\overline{\Omega})$ if $a \in L^{\infty}(\Omega)$,
\item [$(ii)$] $u \in L^{\infty}({\Omega})$ if either $a \in L^{q}(\Omega)\cap L^{\ell^*/(\ell^*+ \alpha - 1)}(\Omega)$ and $0<\alpha \leq1$ or $a \in L^q(\Omega)$ and $\alpha > 1$, where $N/\ell< q \leq q(\alpha)$ with
\begin{equation}\label{probq}
q(s):=\left\{\
\begin{array}{l}
\ell^*/s~\mbox{if}~0<s\leq 1 , \\
(\ell^*+(\alpha-1)\ell^*/\ell)/s~\mbox{if}~s>1,
\end{array}
\right.
\end{equation}
\item [$(iii)$] there exists an only  solution $u \in \w$ of  Problem $(\ref{prob})$   in the sense of $(\ref{final2})$.
\end{enumerate} 
\end{cor}
\nd We are going to take advantage of our techniques to show existence results to the singular-convex problem 
\begin{equation}\label{prob1}
\displaystyle -\mbox{div}(\phi(|\nabla u|)\nabla u) = \frac{a(x)}{u^{\alpha}} + b(x)u^\gamma~\mbox{in}~\Omega, \;\;
u>0~\mbox{in}~\Omega,~u=0~\mbox{on}~\partial \Omega,
\end{equation}
where $ \alpha,  \gamma >0$.         
     
\begin{thm}\label{Teor-prin1}
{\it	Assume  $(\phi_1)-(\phi_3)$ and $0\leq\gamma <\ell -1$. Assume  in addition that  $\displaystyle {a}{d^{-\alpha}}\in L_{\widetilde\Phi}(\Omega)$ and $0\leq b\in L^{\sigma}(\Omega)$ for some $\sigma>\ell/(\ell - \gamma-1)$. Then problem $(\rm \ref{prob1})$ admits a weak  solution $ u \in \w$ such that $u\geq C d \;   \mbox{in} \; \Omega$  for some constant $C>0$. Besides this, $u \in L^{\infty}({\Omega})$ if $b\in L^{\infty}(\Omega)$, and either $a \in L^{q}(\Omega)\cap L^{\ell^*/(\ell^*+ \alpha - 1)}(\Omega)$ with $0<\alpha\leq 1$ or $a \in L^q(\Omega)$ with $\alpha > 1$, where $N/\ell< q \leq q(\alpha+\gamma)$ and $q(s)$ was defined in $(\ref{probq}).$  }
\end{thm}

\begin{rmk} We note that:
\begin{enumerate}
\item [$(a)$]  solutions of both Theorems can be found by variational arguments in some particular cases,
\item [$(b)$] if $\Psi$ is an  N-function such that $\Phi<\Psi << \Phi_*$, then the conditions 
	\begin{equation}
	\displaystyle {a}{d^{-\alpha}}\in L_{\widetilde\Psi}(\Omega) \nonumber
	~\mbox{and}~
	a \in L^{\widetilde{\Phi}}_{loc}(\Omega)\nonumber
	\end{equation}
	\nd could  be used in our results, instead of
	$${a}{d^{-\alpha}}\in L_{\widetilde\Phi}(\Omega)~\mbox{and}~a \in L^{\infty}_{loc}(\Omega),$$
	respectively.  
\end{enumerate}
\end{rmk}

\section{A family of Auxiliary Problems}\label{auxiliary-pb}
In this section, we are going to ``regularize"  problem  (\ref{prob1}) by considering a perturbation by small $\epsilon>0$ of the singular term in (\ref{prob1}). Of course a regularized form of problem  (\ref{prob}) corresponds to $b=0$. Let us consider
\begin{equation}\label{auxprob}
\left\{\
\begin{array}{l}
\displaystyle-\Delta_\Phi u=\frac{a_\epsilon (x)}{(u+\epsilon)^\alpha}+b_\epsilon (x)u^\gamma~\mbox{in}~\Omega\\ 
u>0~\mbox{in}~\Omega,~u=0~\mbox{on}~\partial \Omega
\end{array}
\right.
\end{equation}
for each $\epsilon>0$ given, where the $L^{\infty}(\Omega)$-functions are defined by
$$
a_\epsilon(x)=\min\{a(x),1/\epsilon\},~~ b_\epsilon(x)=\min\{b(x),1/\epsilon\},~ x \in \Omega.
$$ 
\nd Consider the map   $A := A_{\epsilon}:\w\times\w\longrightarrow\mathrm{R}$, defined by
\begin{equation}\label{operadorA}
A(u,\varphi):=\int_\Omega \Big[\phi(|\nabla u|)\nabla u\nabla\varphi dx-\frac{a_\epsilon(x)\varphi}{(|u|+\epsilon)^\alpha}-b_\epsilon(x)(u^+)^\gamma
\varphi \Big] dx,
\end{equation}

\nd Thus, finding a weak solution of \eqref{auxprob} means to find $u \in \w$ such that
\begin{equation}\label{v-eqn}
A(u,\varphi)=0~ \mbox{ for each}~ \varphi\in \w.
\end{equation}

\begin{prop}\label{32}
 	For each $u\in\w$, the functional $A(u,.)$ is linear and continuous. In particular, the  operator $T:= T_{\epsilon}:\w\longrightarrow W^{-1,\widetilde{\Phi}}(\Omega)$ defined by 
	 $$\langle T(u),\varphi\rangle=A(u,\varphi),~u,\varphi\in \w$$
is linear and continuous, and satisfies
\begin{equation}\label{cont-op}
\|T(u)\|_{W^{-1,\widetilde{\Phi}}}\leq 2\|\phi(|\nabla u|)\nabla u\|_{\widetilde \Phi}+\frac C\epsilon\|a_\epsilon\|_{\widetilde \Phi}+C\|b_\epsilon|u|^\gamma\|_{\widetilde \Phi}.
\end{equation}
\end{prop}
\proof: Let $u,\varphi\in\w$. We shall use below the H\"older inequality and the embedding  $\w\hookrightarrow L_\Phi(\Omega)$: 
\begin{eqnarray}\label{A-bem-def}
	|A(u,\varphi)|&\leq&\int_{\Omega} \big[\phi(|\nabla u|)|\nabla u||\nabla \varphi|+\frac{a_\epsilon(x)|\varphi|}{\epsilon^\alpha}+b_\epsilon(x)(u^+)^\gamma|\varphi|~\big] dx\nonumber\\
	&\leq& 2\|\phi(|\nabla u|)\nabla u\|_{\widetilde \Phi}\| \varphi\| +\frac{2}{\epsilon^\alpha} \|a_\epsilon\|_{\widetilde \Phi}\|\varphi\|_\Phi+2\|b_\epsilon|u|^\gamma\|_{\widetilde \Phi}\|\varphi\|_\Phi\nonumber\\
	&\leq&(2\|\phi(|\nabla u|)\nabla u\|_{\widetilde \Phi}+\frac {C}{ \epsilon^\alpha}\|a_\epsilon\|_{\widetilde \Phi}+C\|b_\epsilon|u|^\gamma\|_{\widetilde \Phi})\| \varphi\|.
\end{eqnarray}
It is enough to show that  $\|b_\epsilon|u|^\gamma\|_{\widetilde \Phi}<\infty$. Indeed, by using the embedding  $L_\Phi(\Omega)\hookrightarrow L^\ell(\Omega)$  and  $\gamma\in(0,\ell-1)$ it follows by  Lemma \ref{lema_naru_*} that
\begin{eqnarray}
	\int_\Omega \widetilde\Phi(b_\epsilon(x)|u^\gamma|)dx &\leq& \max\{\|b_\epsilon\|_\infty^{\frac{\ell}{\ell-1}},\|b_\epsilon\|_\infty^{\frac{m}{m-1}}\}\int_\Omega \widetilde\Phi(|u|^\gamma)dx\nonumber\\
	&\leq& C\left(\int_{u\leq1}+\int_{u\geq1}\right)\widetilde{\Phi}(|u|^\gamma) dx\nonumber\\
	&\leq&C\left(|\Omega|+
	\int_{u\geq1}|u|^{\frac{\gamma\ell}{\ell-1}}dx\right)
	\leq C\left(|\Omega|+
	\int_{u\geq1}|u|^\ell dx\right)\nonumber\\
	&\leq&C\left(|\Omega|+
	\int_\Omega|u|^\ell dx\right)\leq C\left(|\Omega|+
	\|u\|^\ell \right)	,
\end{eqnarray}
where  $C=C(b,\Phi,\epsilon) > 0$ is a constant.  So  $A(u,.)$ is linear and continuous.   The claims about $T$  are  now immediate. \hfill\cqd
\vskip.2cm

\nd By proposition \eqref{32} the problem of finding a weak solution of \eqref{auxprob} reduces to find  $u=u_\epsilon\in\w \setminus \{0\}$
such that  $T(u_\epsilon)=0$.

\Section{Applied Generalized Galerkin Method}

\nd In order to find $u=u_\epsilon\in\w \setminus \{0\}$
such that  $T(u_\epsilon)=0$,  we shall employ a Galerkin like method inspired in arguments found in  Browder \cite{Browder}.

\nd We are going to constrain the operator $T$ to finite dimensional subspaces. 
 As a first step  take  a  $\omega\in \w$ such that 
\begin{equation}
\label{w}
a\omega\neq 0~\mbox{and}~  a\omega\in L^1(\Omega),
\end{equation}

\nd Let $F\subset\w$ be a finite dimensional subspace such that $\omega\in F$. Now,  consider the map  $T_F:F\rightarrow F^{\prime}$ given by $T_F = I_F^{\prime}\circ T\circ I_F$, where 
$$
I_F:(F,\|.\|)\longrightarrow(\w,\|.\|),~~ I_F(u)=u
$$
\nd  and let $I_F^{\prime}$ be the adjoint of $I_F$. So, we have that $T_F=T{\big{|}_{F}}$, because 
$$
\langle T_F u,v\rangle = \langle I_F^{\prime}\circ T\circ I_F u,v\rangle=\langle T\circ I_Fu,I_Fv\rangle=\langle T u,v\rangle,~ u,v\in F,
$$
\nd that is, 
\begin{equation}\label{T_K}
	\langle T_F(u),v\rangle:=\int_\Omega \Big[\phi(|\nabla u|)\nabla u\nabla v-\frac{a_\epsilon(x)v}{(|u|+\epsilon)^\alpha}-b_\epsilon(x)(u^+)^\gamma v \Big ]dx,~u,v\in F.
\end{equation}

\nd The result  below, which is a   consequence of the  Brouwer Fixed Point Theorem (see \cite{lions}),   will play  a central  role in solving the finite dimensional equation  $T_F(u) = 0$.
\begin{prop}\label{prop-principal}
	Assume that $S:\mathrm{R}^s\rightarrow\mathrm{R}^s$ is a continuous map such that $( S(\eta),\eta)>0$,~ $|\eta|=r$ for some $r > 0$,  where $(\cdot,\cdot)$ is the usual inner product in $\mathrm{R}^s$ and $|\cdot|$ is its corresponding norm. Then, there is $\eta_0\in  B_r(0)$ such that $S(\eta_0)=0$.
\end{prop}

\begin{prop}\label{S_K-cont}
	 The operator $T_F$ is continuous. 
\end{prop}
\proof: Let  $(u_n)\subseteq F$ be a sequence such that   $u_n\rightarrow u$ in $F$. Since, the operator  $\Delta_{\Phi}:\w\rightarrow W^{-1,\widetilde \Phi}(\Omega)$  given by 
$$
\langle -\Delta_{\Phi}u,v\rangle:=\int_\Omega \phi(|\nabla u|)\nabla u\nabla v dx,~u,v\in\w,
$$
is continuous (see  \cite[Lemma 3.1]{Fukagai}), we  have that    $\Delta_{\Phi}{\big|}_F$ is also continuous. 

\nd To finish our proof, it remains to show that   $T_F-\Delta_{\Phi}{\big|}_F$ is continuous. By applying Lemma \ref{lema_Phi}  and the embedding $L_\Phi(\Omega)\hookrightarrow L^\ell(\Omega)$, it follows,  by eventually passing to a subsequence, that
\begin{itemize}
	\item[(1)] $u_n\rightarrow u$ a.e. in $\Omega$;
	\item[(2)] there is $h\in L^\ell(\Omega)$ such that $|u_n|\leq h$. 
\end{itemize}
Then for each $v\in\w$, 
$$
\frac{a_\epsilon(x)v}{(|u_n|+\epsilon)^\alpha}\longrightarrow\frac{a_\epsilon(x)v}{(|u|+\epsilon)^\alpha},\quad b_\epsilon(x)(u^+_n)^\gamma v \longrightarrow b_\epsilon(x)(u^+)^\gamma v~a.e.~in~ \Omega.
 $$
On the other hand, since $\widetilde{\Phi}$ is increasing, we obtain
\begin{eqnarray}\label{a}
	\widetilde{\Phi}\left(\left|\frac{a_\epsilon(x)}{(|u_n|+\epsilon)^\alpha}-\frac{a_\epsilon(x)}{(|u|+\epsilon)^\alpha}\right|\right)&\leq &	\widetilde{\Phi}\left(\frac{a_\epsilon(x)}{(|u_n|+\epsilon)^\alpha}+\frac{a_\epsilon(x)}{(|u|+\epsilon)^\alpha}\right)\nonumber\\
	&\leq&	\widetilde{\Phi}\left(\frac{2a_\epsilon(x)}{\epsilon^\alpha}\right)\in L^1(\Omega),
\end{eqnarray}
because $0 \leq a_\epsilon \leq 1/\epsilon$. So, 
 by Lebesgue's Theorem, 
$$\int_\Omega\widetilde{\Phi}\left(\left|\frac{a_\epsilon(x)}{(|u_n|+\epsilon)^\alpha}-\frac{a_\epsilon(x)}{(|u|+\epsilon)^\alpha}\right|\right)dx\rightarrow 0,
$$
and as a consequence of $\widetilde{\Phi}\in \Delta_2$, we have 
$$
\left\|\frac{a_\epsilon(x)}{(|u_n|+\epsilon)^\alpha}-\frac{a_\epsilon(x)}{(|u|+\epsilon)^\alpha}\right\|_{\widetilde{\Phi}}\rightarrow 0.
$$

\nd By applying the H\"older inequality, we find that
$$ 
\left|\int_\Omega\left(\frac{a_\epsilon(x)}{(|u_n|+\epsilon)^\alpha}-\frac{a_\epsilon(x)}{(|u|+\epsilon)^\alpha}\right)vdx\right|\leq 2\left\|\frac{a_\epsilon(x)}{(|u_n|+\epsilon)^\alpha}-\frac{a_\epsilon(x)}{(|u|+\epsilon)^\alpha}\right\|_{\widetilde{\Phi}}\|v\|_{{\Phi}}\rightarrow0
$$
\nd for each $v\in\w$. 

\nd Estimating as in  \eqref{a}, we have
\begin{eqnarray}
	\widetilde{\Phi}\left(b_\epsilon|(u_n^+)^\gamma-(u^+)^\gamma|\right)&\leq&\widetilde{\Phi}\left(2|b_\epsilon|_\infty\frac{(u_n^+)^\gamma+(u^+)^\gamma}{2}\right)\nonumber\\	
	&\leq& C\left(\widetilde{\Phi}((u_n^+)^\gamma)+\widetilde{\Phi}((u^+)^\gamma)\right)\nonumber\\
	&\leq&	C\left(|u|^\ell+|h|^\ell+2\right)\in L^1(\Omega),
	\end{eqnarray}
for some $C=C(a,\Phi,\epsilon)>0$. 

\nd Arguing as above, we obtain 
$$
\int_\Omega b_\epsilon(x)[(u_n^+)^\gamma-(u^+)^\gamma]vdx \longrightarrow0
$$
showing  that  $T_F$ is continuous.   \hfill\cqd

\begin{prop}\label{raizT_K}
	There exists  $0\neq u=u_F=u_{\epsilon,F}\in F$ such that
$
 T_F(u) = 0
$	for each $\epsilon>0$ sufficiently small.
\end{prop}

\nd \proof: Let $s:=\dim F$ be the dimension of the subspace $F$, and set   $F=\langle e_1,e_2,...,e_s \rangle$. That is, each  $u\in F$ is uniquely expressed  as 
$$
u=\sum_{j=1}^s\xi_je_j,~ \xi = (\xi_1, \xi_2, \cdots, \xi_s) \in \mathrm{R}^s.
$$
Set $ |\xi| := \|u\|$ and  consider the map  $i=i_F:(\mathrm{R}^s,|.|)\rightarrow(F,\|.\|)$ given by $i(\xi)=u$. 

\nd So,  it follows by Proposition \ref{S_K-cont} and the fact that $i $ is an isometry that the operator   $S_F:\mathrm{R}^s\rightarrow\mathrm{R}^s\nonumber$ given by
\begin{equation}\label{39}
		S_F:=i^{\prime}\circ T_F\circ i
\end{equation} 
is continuous as well, where  $i^{\prime}$ is the adjoint of  $i$.

\nd Besides this, by setting $u := i(\xi)$ for $\xi\in\mathrm{R}^s$, it follows from $(\phi_3)$ and the embeddings $\w\hookrightarrow L_\Phi(\Omega)\hookrightarrow L^\ell(\Omega)\hookrightarrow L^{\gamma+1}(\Omega)$ that
\begin{eqnarray}
\label{r0}
	(S_F\xi,\xi)&=&(i^{\prime}\circ T_F\circ i(\xi),\xi)
	=\langle T_F(u),u\rangle\nonumber\\
	&\geq&\int_\Omega \Big[\phi(|\nabla u|)|\nabla u|^2-\frac{a_\epsilon(x)|u|}{\epsilon^\alpha}-b_\epsilon(x) |u|^{\gamma+1}\Big]dx\nonumber\\
	&\geq&\ell\int_\Omega \Phi(|\nabla u|)dx-\frac{1}{\epsilon^\alpha}\|a_\epsilon\|_{\widetilde{\Phi}}\|u\|_\Phi-|b_\epsilon|_\infty|u|^{\gamma+1}_{\gamma+1}\\
	&\geq&\ell\min\{\|u\|^\ell,\|u\|^m\}-C_1\|u\|-C_2\|u\|^{\gamma+1},\nonumber
\end{eqnarray}
for some positive constants $C_1=C_1(\epsilon)$ and $C_2=C_2(\epsilon)$. So, we can 
choose an $ r_0 =r_0 (\epsilon)> 1$ such that $\ell r_0^\ell-C_1r_0 -C_2r_0^{\gamma+1}>0$. More specifically, for each $\xi$ such that $|\xi| = r_0$, we have
$
	(S_F\xi,\xi)> 0.
$

\nd By the above, it follows from Proposition \ref{prop-principal} that there exists a  $\xi_F\in \overline B_{r_0}(0)$ such that $S_F(\xi_F)=0$, that is, letting  $u=u_F=i(\xi_F)$, it follows from (\ref{39}), that
$$
\langle T_F(u),v\rangle=(S_F(\xi_F),\eta)=0~\mbox{for all}~v\in F,$$
where  $v=i(\eta)$, and hence   $T_F(u)=0$. As a consequence of this, we have
$$
\int_\Omega\Big[\phi(|\nabla u|)\nabla u\nabla v-\frac{a_\epsilon(x)v}{(|u|+\epsilon)^\alpha}-b_\epsilon(x)(u^+)^\gamma v\Big]dx = 0~\mbox{for all}~v\in F.
$$
\nd  We claim that $ u =u_\epsilon\ne 0$ for enough small $\epsilon>0$. Indeed, otherwise by taking $v=w$ and using Lebesgue's Theorem, we obtain
$$
 \int_\Omega {a(x)}wdx = \lim_{\epsilon\to 0}\int_\Omega {a_\epsilon(x)}wdx = 0,
$$
but this is impossible by  (\ref{w}). This ends the proof of Proposition \ref{raizT_K}.
\hfill\cqd
\vskip.2cm

\nd The result below is a direct consequence of the Proposition proved just above.
\begin{cor}\label{r_0}
	 The number $r_0>0$ and the function $u_F \in F$ found  above  satisfy:
$
\|u_F\|\leq r_0,
$	
  $T_F(u_F)=0$, and $r_0>0$  does not depends on subspace $F\subset \w$ with $0<\dim F<\infty$. Besides this, we can choose it independent of $\epsilon>0$ as well if $0<\alpha \leq 1$, $a \in  L^{\ell^*/(\ell^*+ \alpha - 1)}(\Omega)$, and $ b\in L^{\sigma}(\Omega)$ for some $\sigma>\ell/(\ell - \gamma-1)$.
\end{cor}
\proof: The first part of it was proved above. To show that $r_0$ does not depends on $\epsilon>0$, we just redo the estimatives in $(\ref{r0})$
by using  the hypotheses on $a$ and $b$.
\hfill\cqd
\vskip.2cm

\nd Our aim below is  to build a non-zero vector $ u_\epsilon \in \w$ such that $T(u_\epsilon)=0$, where $T$ was given by  Proposition \ref{32}. This will provide us with some  $u_{\epsilon}\in \w$ such that
\begin{equation}
\label{lem31}
\int_\Omega\Big[\phi(|\nabla u|)\nabla u\nabla \varphi-\frac{a_\epsilon(x)\varphi}{(|u|+\epsilon)^\alpha}-b_\epsilon(x)(u^+)^\gamma \varphi\Big]dx = 0,~\varphi\in \w.
\end{equation}
\nd In this direction we have

\begin{lem}\label{raizT}
	There is a non-zero vector $u_{\epsilon} \in \w$ such that  $T(u_{\epsilon}) = 0$ or equivalently $(\ref{lem31})$  holds true.
	\end{lem}

\nd \proof: Let  $w$ as in $(\ref{w})$ and set
$$
\mathcal{A}=\left\{F\subset \w~|~ F~\mbox{is a  finite dimensional subspace of}~ \w~\mbox{and}~ \omega\in F \right\},  
$$

\nd We assume that $\mathcal{A}$ is partially ordered  by  set inclusion.  Take $F_0\in\mathcal{A}$ and  set
$$
V_{F_0}=\left\{u_F\in F~\big|~F\in \mathcal{A},~F_0\subset F,~T_F(u_F)= 0~\mbox{and}~\|u_F\|\leq r_0\right\}.
$$
\nd Note that by Proposition \ref{raizT_K} and Corolary \ref{r_0}, $V_{F_0}\neq \emptyset$.
\vskip.1cm

\nd  Since $V_{F_0} \subset \overline{B_{r_0}}(0)$, then 
$\overline V^\sigma_{F_0} \subset   \overline{B_{r_0}}(0)    $, where  $\overline V^\sigma_{F_0}$ denotes  the weak closure   of $V_{F_0}$. As a matter of this fact,  $\overline V^\sigma_{F_0}$  is weakly compact. \nd Consider the family
 $$
 \mathcal{B}:=\left\{\overline V_{F}^\sigma~|~F\in \mathcal{A}\right\}.
 $$

\nd {\bf Claim.} $\mathcal{B}$ has the finite intersection property.
\vskip.1cm

\nd \nd Indeed, let $\{\overline V_{F_1}^\sigma,\overline V_{F_2}^\sigma,...,\overline V_{F_p}^\sigma \}$ be a finite subfamily of $\mathcal{B}$ and set 
$$
F:=\mbox{span}\{F_1,F_2,...,F_p\}.
$$
\nd By the very definition of $V_{F_i}$, we have that $u_F\in\overline V_{F_i}^\sigma,~i=1,2,...,p$,  that is
$$
\bigcap_{i=1}^p\overline V_{F_i}^\sigma\neq\emptyset.
$$
This ends the proof of the {\bf Claim.}
\vskip.1cm

\nd Since  ${\overline {B}}_{r_0}$ is weakly compact, it follows that (cf. \cite[Thm. 26.9]{munkres})
$$
W:=\bigcap_{F\in\mathcal{A}}\overline V_{F}^\sigma\neq \emptyset.
$$
\nd Let $u_{\epsilon} \in W$. Then $u_{\epsilon} \in \overline V_{F}^\sigma$ for each $F \in \mathcal{A}$. 
\vskip.2cm

\nd Take $F_0 \in \mathcal{A}$ such that $\mbox{ span}\{\omega, u_{\epsilon} \} \subset F_0$. Since  $u_{\epsilon} \in \overline V_{F_0}^\sigma$, it follows by   \cite[Thm. 1.5]{Figueiredo} and the definition of $V_{F_0}$ that there  are sequences $(u_n)=(u_{n,\epsilon})\subset V_{F_0}$ and $(F_n)=(F_{n,\epsilon})\subset \mathcal{A}$ such that $u_n\rightharpoonup u_\epsilon$   in $\w$, $u_n \in F_n$, $\|u_n\|\leq r_0$, $F_0\subset F_n$, and 
\begin{eqnarray}\label{u_n}
	\int_\Omega\phi(|\nabla u_n|)\nabla u_n\nabla v dx=\int_\Omega\left(\frac{a_\epsilon(x)}{(|u_n|+\epsilon)^\alpha}+b_\epsilon(x)(u_n^+)^\gamma\right)vdx
\end{eqnarray}
$\mbox{for each}~v \in F_n$.
\vskip.1cm

\nd Now, by eventually taking subsequences and using $\w \stackrel{comp} \hookrightarrow L_\Phi(\Omega)$,  we obtain that
$u_n\rightarrow u_\epsilon$ in $L_\Phi(\Omega)$,~ $u_n\rightarrow u_\epsilon$ a.e. in $\Omega$ and $(|u_n|)$ is bounded away by some function in $L_\Phi(\Omega)$.
\vskip.2cm

\nd  Set  $v_n=u_n-u_\epsilon$ and note that $v_n \in F_n$, because $u_n \in F_n$ and $u_\epsilon \in F_0 \subset F_n$ in \eqref{u_n}.
 Then
 \begin{eqnarray}\label{S+}
	\lim\langle -\Delta_{\Phi} (u_n),u_n-u_\epsilon\rangle&=&\lim\int_\Omega \left(\frac{a_\epsilon(x)}{(|u_n|+\epsilon)^\alpha}+b_\epsilon(x)(u_n^+)^\gamma\right)(u_n-u_\epsilon)dx\nonumber\\
	&\leq& \lim\int_\Omega \left(\frac{a_\epsilon(x)}{\epsilon^\alpha}+b_\epsilon(x)|u_n|^\gamma\right)|u_n-u_\epsilon|dx.
\end{eqnarray}

\nd As  $\w\stackrel{comp}\hookrightarrow L_{\Phi}(\Omega)$, we have
$$
\left|\int_{\Omega}\frac{a_\epsilon(x)}{\epsilon^\alpha}(u_n-u_0)dx\right|\leq\frac{1}{\epsilon^\alpha}\|a_\epsilon\|_{\widetilde\Phi}\|u_n-u_\epsilon\|_\Phi\rightarrow0.
$$

\nd Recalling that  $\gamma<\ell-1$,  $\w\hookrightarrow L^\ell(\Omega)$ and $(u_n)$ is bounded in  $L^\ell(\Omega)$, we get
\begin{eqnarray*}
	\int_\Omega b_\epsilon(x)|u_n|^\gamma|u_n-u_\epsilon|dx &\leq& |b_\epsilon|_\infty\left(\int_\Omega|u_n|^{\frac{\gamma\ell}{\ell-1}} dx\right)^{\frac{\ell-1}{\ell}} |u_n-u_\epsilon|_\ell  \nonumber\\
	&\leq&|b_\epsilon|_\infty\left(|\Omega|+\int_\Omega|u_n|^\ell dx\right)^{\frac{\ell-1}{\ell}}|u_n-u_\epsilon|_\ell\rightarrow 0.\nonumber\\
\end{eqnarray*}

\nd Now, by using the facts above, it follows from (\ref{u_n}) that
$$
\lim\langle -\Delta_{\Phi}(u_n),u_n-u_\epsilon\rangle\leq 0,
$$
and a consequence of this, we have that $u_n\rightarrow u_\epsilon$ in $\w$, because $-\Delta_{\Phi}$ satisfies the condition  $(S_+)$ (see \cite[Prop. A.2]{JVMLED}). 

\nd So, passing to a subsequence if necessary, we have
\begin{itemize}
	\item[(1)] $\nabla u_n\rightarrow \nabla u_\epsilon$ a.e. in $\Omega$,	
	\item[(2)] there is  $h\in L_\Phi(\Omega)$ such that $|\nabla u_n|\leq h$.
\end{itemize}
\nd Since $\varphi \in \w$, it follows of the fact that $t \phi(t)$ is nondecreasing in $[0,\infty)$ and (2), that
\begin{eqnarray}
	|\phi(|\nabla u_n|)\nabla u_n\nabla \varphi |&\leq& \phi(|\nabla u_n|)|\nabla u_n||\nabla \varphi|
	\leq \phi(h)h|\nabla \varphi|\nonumber\\
	&\leq& \widetilde{\Phi}(\phi(h)h)+\Phi(|\nabla \varphi|)\leq\Phi(2h)+\Phi(|\nabla \varphi|)\in L^1(\Omega),\nonumber
\end{eqnarray}
that is, it follows by the Lebesgue Theorem, that
$$
\int_\Omega\phi(|\nabla u_n|)\nabla u_n\nabla \varphi dx \longrightarrow\int_\Omega\phi(|\nabla u_\epsilon|)\nabla u_\epsilon\nabla \varphi dx.
$$	
\vskip.1cm

\nd Now, by passing to the limit in \eqref{u_n} and using the above informations, we get that $u_\epsilon$ satisfies (\ref{lem31}), that is, $T_\epsilon(u_\epsilon)=T(u_\epsilon)=0$ for each $\epsilon>0$, since   $\varphi \in \w$ was taken arbitrarily.  By arguments as in the proof of Proposition \ref{raizT_K} we infer that 
 $u_\epsilon\not\equiv 0$. \hfill\cqd

\begin{lem}\label{solu-aux}
 The function $u_\epsilon\in C^{1,\alpha_\epsilon}(\overline \Omega)$, for some $0<\alpha_\epsilon\leq 1$, and it is a solution of \eqref{auxprob}.
\end{lem}


\nd \proof: By Lemma \ref{raizT}, it remains  to show that  $u_\epsilon> 0$. Set $-u^-_\epsilon$ as a test function in  \eqref{lem31}. So, it follows by Remark \ref{rem51} (see Appendix), that
\begin{eqnarray}
		\ell\int_\Omega\Phi(|\nabla u^-_\epsilon|)dx&\leq& \int_\Omega\phi(|\nabla u^-_\epsilon|)|\nabla u^-_\epsilon|^2dx\nonumber\\
		&=&-\int_\Omega \frac{a_\epsilon(x)}{(|u^-_\epsilon|+ \epsilon)^\alpha}u^-_\epsilon dx,
\end{eqnarray}

\nd which implies that $u^-_\epsilon\equiv0$. So, $u_\epsilon$ satisfies
\begin{eqnarray}\label{solfrac}
\int_\Omega\phi(|\nabla u_\epsilon|)\nabla u_\epsilon\nabla\varphi \!=\!\int_\Omega\frac{a_\epsilon(x)}{(u_\epsilon+\epsilon)^\alpha}\varphi dx\!
+\!\int_\Omega b_\epsilon(x)u_\epsilon^\gamma\varphi dx, \varphi\in \w.
\end{eqnarray}

\nd Finally, for  each $p\in(m,\ell^*)$, it follows that
$$
|f(x,t)|:=\frac{a_\epsilon(x)}{(|t|+\epsilon)^\alpha}+b_\epsilon(x)(t^+)^\gamma\leq C_\epsilon (1+|t|^{p-1})~\mbox{and}~ \displaystyle\lim_{t\rightarrow\infty} \frac{t^p}{\Phi_*(\lambda t)}=0
$$
\nd for each $\epsilon > 0$ given. So by  \cite[Corollary 3.1]{fang}, $u_\epsilon\in C^{1,\alpha_\epsilon}(\overline \Omega)$ for some $0<\alpha_\epsilon\leq 1$. Now, by summing up the term $u_\epsilon\phi(u_\epsilon)$ to both sides of  
(\ref{solfrac}) and applying  \cite[Proposition 5.2]{JVMLED} we infer that $u_\epsilon>0$.  In conclusion, $u_\epsilon$ is a solution of \eqref{auxprob}.  \hfill\cqd

\Section{Comparison of Solutions and Estimates}

\nd Let $n \geq 1$ be an integer and take $\epsilon=1/n$. Let $u_n\in \w\cap C^{1,\alpha_n}(\overline{\Omega})$, for some $\alpha_n \in (0, 1]$, denotes the solution of (\ref{auxprob}), both for $b=0$ and $b\geq 0$ not identically null, given by Lemma \ref{solu-aux}, that is,  
\begin{equation}\label{ineq-2}
\displaystyle-\Delta_\Phi u_n =\frac{a_n(x)}{(u_n+1/n)^\alpha} + b_n(x) u^{\gamma}~~\mbox{in}~~\Omega,~~ 
u_n>0~\mbox{in}~\Omega,~u_n=0~\mbox{on}~\partial \Omega
\end{equation}
\vskip.1cm

\nd We have the following result on comparison of solutions.

\begin{lem}\label{est}
	The following inequalities hold:
\vskip.2cm

\nd $\mbox{\rm (i)}$~ $u_n+1/n\geq u_1$ for each integer $n \geq 1$
\vskip.2cm	
	
\nd $\mbox{\rm (ii)}$~ $u_1\geq C d~\mbox{in}~ \Omega$~  for some $C > 0$ which independs of $n$.
\end{lem}
\vskip.3cm

\nd \proof. First, we consider $b=0$ in (\ref{ineq-2}), that is,
\begin{equation}\label{ineq-1}
\displaystyle-\Delta_\Phi u_n =\frac{a_n(x)}{(u_n+1/n)^\alpha} ~~\mbox{in}~~\Omega,~~ 
u_n>0~\mbox{in}~\Omega,~u_n=0~\mbox{on}~\partial \Omega.
\end{equation}
 So, by \eqref{ineq-1} we have
\begin{equation}\label{serrin-1}
\displaystyle \mbox{div}(\phi(|\nabla u_1|)\nabla u_1) - \frac{a_1(x)}{(u_1+1)^\alpha} \geq  0~~\mbox{in}~~\Omega, 
\end{equation}
in the weak sense.
On the other hand, since  

\begin{equation*}
 \frac{a_n(x)}{(w_n+1/n)^\alpha} \geq \frac{a_1(x)}{( (w_n+1/n) + 1)^\alpha }~\mbox{in}~\Omega.
\end{equation*}

\nd  we get by \eqref{ineq-1} that 
\begin{equation}\label{serrin-2}
\displaystyle \mbox{div}(\phi(|\nabla  (u_n + 1/n)|)\nabla  (u_n + 1/n)) - \frac{a_1(x)}{( (u_n + 1/n)+1)^\alpha} \leq  0~~\mbox{in}~~\Omega, 
\end{equation}
\nd in the weak sense, (test finctions are taken non-negative).

\nd By applying  Theorem 2.4.1 in  \cite{serrin-pucci} to (\ref{serrin-1}) and  (\ref{serrin-2}), we obtain  $u_n+1/n\geq u_1$.
\vskip.1cm
 
\nd Now,  since $\partial \Omega$  is smooth, it follows  by \cite[Lemma 14.16]{Gilbarg} that the distance function $x \mapsto d(x)$  satisfies
$$
d \in C^2(\overline{\Omega}),~ d > 0~\mbox{on}~\overline{\Omega}_\delta~\mbox{and}~  \frac{\partial d}{\partial \eta}<0~\mbox{on}~\overline{\Omega} \setminus \Omega_\delta,
$$ 
where $\Omega_\delta=\{x\in\overline\Omega~|~d(x)>\delta\}$ for some $\delta > 0$,  and $\eta$ stands for the  exterior unit normal to  $\partial \Omega$. 
\vskip.1cm

\nd Now, since $u_1 \in \w\cap C^{1,\alpha_1}(\overline{\Omega})$ is a solution of

\begin{equation}\label{aux-01} 
-\Delta_\Phi u=\displaystyle\frac{a_1(x)}{(u+1)^\alpha}~\mbox{in}~\Omega,~~ 
u>0~\mbox{in}~\Omega,~u=0~\mbox{on}~\partial \Omega,
\end{equation}

\nd   it follows by  \cite[Lemma 4.2]{fang} that 
$$
\frac{\partial u_1}{\partial \eta}<0~\mbox{on}~\overline{\Omega} \setminus \Omega_\delta.
$$
\nd  So there is  a constant $C>0$ such that 
$$\frac{\partial u_1}{\partial \eta} \leq C\frac{\partial d}{\partial \eta}~\mbox{on}~\overline{\Omega} \setminus \Omega_\delta,$$
\nd and as a consequence
\begin{eqnarray}\label{d(x)} 
C d(x) \leq u_1(x)~\mbox{for}~ x \in \Omega.
\end{eqnarray}

\nd This ends the proof of Lemma \ref{est} for $b=0$. If $b$ is not identically null, we redo the above proof by considering (\ref{serrin-1}) and obtaining (\ref{serrin-2}) as a consequence of $b$ be non-negative. This ends the proof.\hfill\cqd
\vskip.3cm

\nd We have the following estimates.

\begin{lem}\label{est-2}
	\label{lem-thm-2.1}
Let $u_n \in C^{1,\alpha_n}(\overline{\Omega}) $ be a solution of $(\ref{ineq-1})$. Then there is a constant $C > 0$ such that 
	$$
	\Vert[(u_n+1/n)^{({\alpha+\ell-1)}/{\ell}} - (1/n)^{{(\alpha+\ell-1)}/{\ell}} ] \Vert_{1,\ell} \leq C,~\mbox{for all integer}~n \geq 1.
	$$
\end{lem}
\nd where $\Vert \cdot \Vert_{1,\ell}$  above is the norm of $W_0^{1,\ell}$.
\vskip.2cm

\nd \proof.  At first notice that  
$$
u_n,~ [(u_n + 1/n)^{\alpha}-(1/n)^{\alpha}] \in \w\cap C^{1,\alpha_n}(\overline{\Omega}) \subset W_0^{1,\ell}(\Omega).
$$ 
\nd By estimating, we get
$$
\displaystyle \ell \alpha \Phi(1)\int_{|\nabla u_n|\geq 1} |\nabla u_n|^\ell(u_n+\frac 1n)^{\alpha-1} dx \leq
$$
$$
\displaystyle   \ell\alpha \Phi(1) \big[\int_{|\nabla u_n|< 1}|\nabla u_n|^m (u_n+\frac 1n)^{\alpha-1} dx
+\int_{|\nabla u_n|\geq 1} |\nabla u_n|^\ell(u_n+\frac 1n)^{\alpha-1} dx \big]  \leq 
$$
$$
\displaystyle \ell\alpha \Phi(1) \int_\Omega \min\{|\nabla u_n|^\ell,|\nabla u_n|^m\}(u_n+ 1/n)^{\alpha-1} dx.
$$
\nd Applying Remark $\ref{rem51}$ and Lemma \ref{lema_naru}, (both in the Appendix), and using $[(u_n+1/n)^{\alpha} - (1/n)^{\alpha}]$ as a test function in 
\eqref{ineq-1}, we find
\begin{eqnarray} \nonumber
\!\!\!\!\!\!\ell\alpha \Phi(1)\int_{|\nabla u_n|\geq 1}|\nabla u_n|^\ell (u_n\!\!\!&+&\!\!\!1/n)^{\alpha-1}dx \leq     \ell\alpha\int_{\Omega}\Phi(|\nabla u_n|)(u_n+1/n)^{\alpha-1}dx  \\ \nonumber 
&\leq& \alpha\int_{\Omega}\phi(|\nabla u_n|)|\nabla u_n|^2(u_n+1/n)^{\alpha-1}dx\\ \nonumber
&=& \int_{\Omega} \frac{a_n(x)[(u_n+1/n)^\alpha- (1/n)^\alpha]}{(u_n+1/n)^\alpha}dx \\ 
&\leq& |a|_1. \label{461}
\end{eqnarray}
\vskip.2cm

\nd When $\alpha \leq 1$, it follows from Lemma $\ref{est}$, that
\begin{eqnarray}\label{473}
\!\!\!\!\!\!\ell\alpha \Phi(1)\int_{|\nabla u_n|\leq 1}\!\!\!\!\!\!|\nabla u_n|^\ell (u_n+1/n)^{\alpha-1}dx  \leq \\ \nonumber
\ell\alpha \Phi(1)\Big[ |\Omega|+ C^{\alpha-1}\!\!\!\int_{\Omega}d(x)^{\alpha-1} \Big]:=D
\end{eqnarray}

\noindent  which is finite, by a well known result, cf.  Lazer and McKenna $\cite{lazerMcKenna-1}$.
\vskip.2cm

\nd It follows from  $(\ref{461})$ and $(\ref{473})$, that
\begin{eqnarray}
\int_{\Omega}\left|\nabla \left((u_n+1/n)^{\frac{\alpha-1+\ell}{\ell}}\right)\right|^\ell dx 
&\leq& \left(\frac{\alpha+\ell-1}{\ell}\right)^\ell \frac{1}{\ell\alpha \Phi(1)}(\Vert a \Vert_1+ D),\nonumber
\end{eqnarray}
because
$$\left|\nabla \left((u_n+1/n)^{\frac{\alpha+\ell-1}{\ell}}\right)\right|^\ell=\left(\frac{\alpha+\ell-1}{\ell}\right)^\ell |\nabla u_n|^\ell (u_n+1/n)^{\alpha-1}.$$
Hence, $[(u_n+1/n)^{({\alpha+\ell-1)}/{\ell}} - (1/n)^{{(\alpha+\ell-1)}/{\ell}}]$ is bounded in $W_0^{1, \ell}(\Omega)$.  
\vskip.3cm

\noindent  When $\alpha > 1$,  we have
\begin{eqnarray}
\label{472}
\!\!\!\!\!\!\ell\alpha \Phi(1)\int_{|\nabla u_n|\leq 1}|\nabla u_n|^\ell (u_n+1/n)^{\alpha-1}dx  \leq \\ \nonumber \ell\alpha \Phi(1)\Big[ |\Omega|+\int_{u_n> 1} (u_n+1/n)^{\alpha-1}dx. \Big]
\end{eqnarray}

\nd Summing up $(\ref{461})$ and $(\ref{472})$, we obtain a positive constant $C$ such that
\begin{eqnarray}\label{est-u_n}
\int_{\Omega}\left|\nabla \left((u_n+1/n)^{\frac{\alpha-1+\ell}{\ell}}\right)\right|^\ell dx\leq C\left(1+\int_{u_n>1}(u_n+1/n)^{\alpha-1}dx\right).
\end{eqnarray}

\nd Now, by picking $\epsilon$ such that   $0<\epsilon < \ell-{\ell(\alpha-1)}/({\alpha+\ell-1})$, it follows from \eqref{est-u_n}, using $u_n>1$ and of the embbeding  $W_0^{1,\ell}(\Omega)\hookrightarrow L^\ell(\Omega)\hookrightarrow L^{\ell-\epsilon}(\Omega)$, that
\begin{eqnarray}
\left\|\nabla \left((u_n+1/n)^{\frac{\alpha-1+\ell}{\ell}}\right)\right\|^\ell_\ell &\leq& \nonumber
  C\left(1+\int_{u_n>1}\left((u_n+1/n)^{\frac{\alpha+\ell-1}{\ell}}\right)^{\ell-\epsilon}dx\right) \\ \nonumber
&\leq&   C\left(1+\left\|\nabla\left((u_n+1/n)^{\frac{\alpha+\ell-1}{\ell}}\right)\right\|^{\ell-\epsilon}_\ell\right),\nonumber
\end{eqnarray}
for  some $C>0$. That is, $[(u_n+1/n)^{({\alpha+\ell-1)}/{\ell}} - (1/n)^{{(\alpha+\ell-1)}/{\ell}}]$ is bounded in $W^{1,\ell}_0(\Omega)$ as well. This ends the proof of Lemma \ref{est-2}.\hfill\cqd

\section{Proof of The Main Results}

We begin proving Theorem \ref{Teor-prin} that treats about existence of positive solution to the pure singular problem \ref{prob}.

\subsection{Pure Singular Problem - Existence of Solutions}
\nd {\bf Proof of {\rm (i)} of Theorem $\ref{Teor-prin}$}~~Assume first that $ad^{-\alpha} \in L_{\tilde{\Phi}}(\Omega)$. 
Since $u_n \in \w$ satisfies $(\ref{ineq-1})$, it follows from Remark $\ref{rem51}$, Lemma $\ref{lema_naru}$,  $\eqref{d(x)}$ and H\"older inequality, that
\begin{eqnarray}\label{a-tende0}
\ell \zeta_0(\Vert \nabla u_n \Vert_{\Phi}) &\leq& \ell \int_{\Omega} \Phi(\vert \nabla u_n \vert) dx \leq \int_{\Omega} \phi(\vert \nabla u_n \vert) \vert \nabla u_n \vert^2  dx \nonumber\\
	&=&\int_{\Omega}\frac{a_n(x)}{(u_n+\frac{1}{n})^\alpha}u_ndx\leq C\int_{\Omega}\frac{a(x)}{d^\alpha}|u_n|dx\nonumber\\
	&=&C\left(\int_{\Omega/\Omega_\delta}+\int_{\Omega_\delta}\right)\frac{a(x)}{d^\alpha}|u_n|dx\\
	&\leq& C \int_{\Omega}|u_n|dx+C\int_{\Omega}\frac{a(x)}{d^\alpha(x)}|u_n|dx\nonumber\\
	&\leq& C \|u_n\|_\Phi+2C\left\|\frac{a}{d^\alpha}\right\|_{\widetilde{\Phi}}\|u_n\|_\Phi,\nonumber
	\end{eqnarray}
	where we used $a_n \leq a$ just above. It follows from our assumptions and from $\w \stackrel{\tiny{cpt}} \hookrightarrow L_\Phi(\Omega)$, that $(u_n)\subset \w$ is bounded. If $0<\alpha\leq1$ and $a \in  L^{\ell^*/(\ell^*+ \alpha - 1)}(\Omega)$, then the boundedness of $(u_n)$ in $\w$ is a consequence of Corollary \ref{r_0}.
	
So, in both cases, up to subsequences, there exist  $u\in \w$ and $\theta\in L_\Phi(\Omega)$ such that
$$
	(1)~  u_n\rightharpoonup u~ \mbox{in}~  \w,~~
	(2)~ u_n\rightarrow u~ \mbox{in}~ L_\Phi(\Omega),~~(3)~ u_n\rightarrow u~ \mbox{a.e. in}~\Omega,
$$
$$
(4)~ 0\leq u_n\leq \theta.
$$
As a first consequence of these facts, it follows from Lemma $\ref{est}$ and $(3)$  that $u\geq Cd~ \mbox{a.e. in}~ \Omega$.

\nd Now, by  using $u_n - u$ as a test function in \eqref{ineq-1} and following similar arguments as in $(\ref{S+})$, we get 
\begin{eqnarray}\label{a-tende0}
\langle -\Delta_{\Phi} u_n,u_n-u\rangle&\leq &	\left|\int_{\Omega}\frac{a_n(x)}{(u_n+1/n)^\alpha}(u_n-u)dx\right|\ \nonumber\\
	&\leq& \Big[C+2\left\|\frac{a}{d^\alpha}\right\|_{\widetilde{\Phi}}\Big]\|u_n-u\|_\Phi
	\end{eqnarray}
for some $C>0$ independent of $n$.  Since, the operator  $-\Delta_{\Phi}$ is of the type $S_{+}$, it follows from $(2)$ and $(\ref{a-tende0})$ that $u_n\rightarrow u~\mbox{in}~\w$.  
\vskip.2cm

\nd To finish our proof, given  $\varphi \in \w$, it follows from Lemma $\ref{est}$, that
\begin{eqnarray}
\Big|\frac{a_n}{(u_n+1/n)^\alpha}\varphi \Big| \leq \frac{a}{d^{\alpha}}\Big(\frac{d}{u_n+1/n}\Big)^{\alpha} {|\varphi|}\leq C\frac{a}{d^{\alpha}} {|\varphi|}\in L^1(\Omega),\nonumber
\end{eqnarray}
that is, by passing to the limit in~\eqref{ineq-1}, we obtain that $u$ is a solution of \eqref{prob}. This ends our proof.
  \hfill\cqd
\vskip.5cm

\nd We were not able to employ the above arguments in the proof of {\rm (ii)} of Theorem \ref{Teor-prin}, because in such case  we do not know if $a/d^{\alpha}$  belongs to $L_{\widetilde{\Phi}}(\Omega)$, that is, the sequence $(u_n)$ likely is not bounded in $\w$. Instead,  it was possible to show  that $(u_n)$ is bounded in $W_{loc}^{1,\Phi}(\Omega)$. This was done by applying Lemma \ref{est-2}. 
\vskip.3cm

\nd {\bf Proof of {\rm (ii)} of Theorem $\ref{Teor-prin}$.} Given $U\subset\subset \Omega$, let $\delta_U = \min\{d(x)~/~x \in U\}>0$. So, it follows from Lemma \ref{est}, that 
$$u_n +1/n \geq C \delta_U:= C_U>0~\mbox{in}~U,$$
that is, for $n >1$ enough big, we can take $(u_n +1/n -C_U)^+$ as a test function in (\ref{ineq-1}), to obtain

\begin{eqnarray}
\label{46}
\begin{array}{lll}
\displaystyle \int_{U} \phi(\vert \nabla u_n \vert)
\vert \nabla u_n \vert^2 &\leq& \displaystyle \int_{u_n+1/n\geq C_U} \phi(\vert \nabla u_n \vert)
\vert \nabla u_n \vert^2 dx\\
&\leq & \displaystyle\int_{u_n+1/n \geq C_U} \frac{a(x)}{(u_n + 1/n)^{\alpha-1}} dx\\
 &\leq &\displaystyle  \frac{1}{C_U^{\alpha -1}} \int_{ \Omega}a dx <\infty,
 \end{array}
\end{eqnarray}
 because $a \in L^1(\Omega)$, and $\alpha \geq 1$.  

So, it follows from  Remark \ref{rem51} and Lemma \ref{lema_naru},  that $(u_n) \subset W^{1,\Phi}(U)$ is bounded.  That is,   there exist $(u^U_{n_{1}}),u^U \in W^{1,\Phi}(U)$  such that
 $u^U_{n_{1}}\rightharpoonup u^U$ in $W^{1,\Phi}(U)$,
 $u^U_{n_{1}}\to u^U$ in $L_{\Phi}(U)$, 
 $u^U_{n_{1}}(x) \to u^U(x)$ a.e. in $U$. In particular, it follows from  Lemma \ref{est} and of the pointwise convergence that $u \geq Cd~ \mbox{a.e. in}~ U$.
 
 Hence, by using a  Cantor diagonalization argument applied
to an exhaustion ${U_k}$ of $\Omega$ with $U_k\subset \subset U_{k+1} \subset \subset \Omega$, we show that there is $u \in W^{1,\Phi}_{loc}(\Omega)$ such that $u_k \to u$ in $W^{1,\Phi}_{loc}(\Omega)$ and $u \geq Cd~ \mbox{a.e. in}~ \Omega$.
 
 Now, we are going to show that this $u$ satisfies the equation in (\ref{prob}). Given  $\varphi\in C^{\infty}_0(\Omega)$, let $\Theta\subset\subset\Omega$ be the support of $\varphi$. So, by very above informations, we have that
\vskip.2cm
 
$(a)$ $u_n\rightharpoonup u$ in $W^{1,\Phi}(\Theta)$,~~$(b)$ $u_n\to u$ in $L_{\Phi}(\Theta)$,~~ $(c)$ $u_n(x) \to u(x)~ \mbox{a.e. in}~ \Theta$ 
\vskip.2cm

\nd and there exists  $\theta \in L_{\Phi}(\Theta)$ such that $u_n \leq \theta$ in $\Theta$.
\vskip.2cm
 
So, by using $\varphi(u_n-u)$ as a test function in $(\ref{ineq-1})$, $L_{\Phi}(\Theta) \hookrightarrow L^1(\Theta)$, and $(b)$ above, we obtain
\begin{eqnarray}
	\Big\vert \int_\Theta\phi(|\nabla u_n|)\nabla u_n\nabla (\varphi(u_n-u))\Big\vert dx&\leq&\frac{1}{c^{\alpha}_{d}}\int_\Theta a_n \vert \varphi( u_n-u)\vert dx\\\nonumber
	&\leq& C_\varphi \Vert a \Vert_{L_{\widetilde{\Phi}(\Theta)}}\Vert u_n-u \Vert_{L_{\Phi}(\Theta)}\longrightarrow 0,\nonumber
\end{eqnarray}
where $\Theta\subset\subset\Omega$ is the support of $\varphi$. That is, 
\begin{eqnarray}
\label{412}
	\int_\Theta\phi(|\nabla u_n|)\nabla u_n\nabla (u_n-u)\varphi dx=\int_\Theta\phi(|\nabla u_n|)\nabla u_n\nabla \varphi(u_n-u)dx+o_n(1).
\end{eqnarray}
Besides this, it follows from Holder's inequality,  $(b)$ above and the property $\widetilde{\Phi}(\phi(t)t) \leq \Phi(2t)$ for $t>0$, that
\begin{eqnarray}
\left|\int_\Theta\phi(|\nabla u_n|)\nabla u_n\nabla \varphi(u_n-u)\right|dx&\leq&C_\varphi\int_\Theta\phi(|\nabla u_n|)|\nabla u_n||u_n-u| dx\nonumber\\
&\leq&C_\varphi\|\phi(|\nabla u_n|)|\nabla u_n|\|_{L_{\widetilde{\Phi}}(\Theta)}\|u_n-u\|_{L_\Phi(\Theta)}\rightarrow 0 \nonumber\\
&\leq&C_\varphi\|u_n-u\|_{L_\Phi(\Theta)}\rightarrow 0,\nonumber
\end{eqnarray}
and using this information in $(\ref{412})$, we obtain that
\begin{eqnarray}
\label{413}
	\int_\Omega\phi(|\nabla u_n|)\nabla u_n\nabla (u_n-u)\varphi dx=o_n(1).
\end{eqnarray}

\nd Besides this, we note that
\begin{eqnarray}
\left|\int_\Theta\phi(|\nabla u|)\nabla u\nabla (u_n-u) \varphi dx\right|&\leq&\left|\int_\Theta\phi(|\nabla u|)\nabla u\nabla [\varphi(u_n-u) ]\varphi dx \right | \nonumber\\
&+&\left|\int_\Theta\phi(|\nabla u|)\nabla u\nabla \varphi (u_n-u)  dx\right|,\nonumber
\end{eqnarray}
and the first integral on the right side goes to zero, due to $(a)$ above, and the second one converges to zero due to $(b)$ above. That is,
\begin{eqnarray}
\label{415}
\left|\int_\Theta\phi(|\nabla u|)\nabla u\nabla (u_n-u) \varphi dx\right|&\to 0.
\end{eqnarray}

\nd So, it follows from $(4.13)$ and $(4.15)$, that
\begin{eqnarray}
\int_\Theta\Big(\phi(|\nabla u_n|)\nabla u_n-\phi(|\nabla u|)\nabla u,\nabla u_n-\nabla u\Big) \varphi dx \to 0,
\end{eqnarray}
and a consequence of this togheter with the Lemma 6 in $\cite{murat}$, we have that $\nabla u_n(x) \to \nabla u(x)~ \mbox{a.e. in}~ \Theta,$ that is, 
$$
\phi(|\nabla u_n(x)|) \nabla u_n(x) \to \phi(|\nabla u(x)|) \nabla u(x) ~ \mbox{a.e. in}~ \Theta.
$$
\nd In addition, since $(\phi(|\nabla u_n|) \nabla u_n )\subset (L_{\widetilde{\Phi}}(\Theta))^N$ is bounded, it follows from Lemma 2 in $\cite{gossez-Czech}$ that 
$$
\phi(|\nabla u_n|)\nabla u_n \rightharpoonup \phi(|\nabla u|) \nabla u~\mbox{in}~  (W^{1,\Phi}(\Theta))^N.
$$
Now, passing to limit in $(\ref{ineq-1})$, we obtain that $u\in W_{loc}^{1,\Phi}(\Omega)$ satisfies
$$\int_\Omega \phi(|\nabla u|)\nabla u\nabla \varphi dx=\int_\Omega \frac{a(x)}{u^\alpha}\varphi dx.$$
Besides this, it follows from Lemma \ref{est-2}, that  
$$u_n^{\frac{\alpha-1-\ell}{\ell}} \rightharpoonup v~\mbox{in}~ W_0^{1,\ell}(\Omega),$$ that is, $u^{({\alpha-1-\ell})/{\ell}} \in W_0^{1,\ell}(\Omega)$ as well.  This ends our proof. \hfill\cqd
\medskip

Below, we take advantage of the former arguments to show existence of solutions to Problem \ref{prob1}. The greatest effort is done to show $L^{\infty}$-regularity of its solutions. 

\subsection{Convex Singular Problem - Regularity of Solutions}

\nd{\bf Proof of Theorem \ref{Teor-prin1}:}  Since $0 < \gamma <\ell - 1$ and $0\leq a\in L^{q}(\Omega)$ for some $q>\ell/(\ell - \gamma-1)$, it follows by arguments similar to those used in the proof of  Theorem $\ref{Teor-prin}$ that there exist both a sequence of aproximating solutions still denoted by $(u_n)$ and  a corresponding solution $ u \in \w$ to problem $(\ref{prob1})$  such that $u\geq C d \;   \mbox{in} \; \Omega$  for some constant $C>0$.
\vskip.1cm

\nd {\bf Claim.}  $u \in L^{\infty}(\Omega)$.  
\vskip.1cm

\nd The proof of this  {\bf Claim} uses arguments driven by  a Moser Iteration Scheme. Parts of our argument were motivated by reading \cite{GuoGao}. However our proof in the present paper is selfcontained.  In order to show the {\bf Claim}, set
$$
\beta_1:=(\ell+\alpha-1)q'>0,~~\beta_k^*:=\beta_k+\beta_1,~\beta_{k+1}:=\frac{\ell^*}{\ell q'}\beta_k^*,~~\delta:={\ell^*}/({q'\ell}),
$$
where $1/q' + 1/q=1$.

\nd We point out that $\delta>1$ because $q>{N}/{\ell}$. In addition,
\begin{eqnarray}
\label{beta1}
\beta_k^*=\left(2\delta^{k-1}+\delta^{k-2}+...+1\right)\beta_1=\frac{2\delta^k-\delta^{k-1}-1}{\delta-1}\beta_1,
\end{eqnarray}
\begin{eqnarray}
\label{beta2}
\beta_k=\frac{2\delta^k-\delta^{k-1}-\delta}{\delta-1}\beta_1,
\end{eqnarray}
\nd and since $\delta > 1$,~  $\beta_k \nearrow \infty$. 
\vskip.1cm

\nd Now, taking  $k_0$ such that $\beta_{k_0},\beta_{k_0}+q'(\alpha-1) > 1$, we have that $u_n^{{\beta_k}/({q'}+\alpha)}$ is a test function for each $k\geq k_0$ and using it in  
$\eqref{lem31}$, we obtain
\begin{eqnarray}
\label{moser}
	\frac{\beta_k}{q'}\int_\Omega\phi(|\nabla u_n|)|\nabla u_n|^2u_n^{\frac{\beta_k}{q'}+\alpha-1} dx&\leq&\int_\Omega\Big(\frac{a_n u_n^{\frac{\beta_k}{q'}+\alpha}}{(u_n+1/n)^\alpha }+b u_n^{\frac{\beta_k}{q'}+\alpha+\gamma}\Big)dx\nonumber\\
	&\leq& \int_{\Omega}\Big( au_n^{\frac{\beta_k}{q'}}+b u_n^{\frac{\beta_k}{q'}+\alpha+\gamma}\Big)dx\\
	&\leq& \|a\|_q\|u_n\|_{\beta_k}^{\frac{\beta_k}{q'}}+\|b\|_\infty\|u_n\|_{\beta_k}^{\frac{\beta_k}{q'}}\|u_n^{\alpha+\gamma}\|_{q}~ .\nonumber
\end{eqnarray}

\nd We claim that $\|u_n^{\alpha+\gamma}\|_{q}$ is bounded. Indeed, \vskip.2cm

\nd if  $(\alpha+\gamma)q\leq 1$, it follows that $\alpha\leq 1$, because $q>N/\ell>1$. In this case, it follows from Corollary $ \ref{r_0}$ that  $u_n$ is bounded in $W_0^{1,\Phi}(\Omega)$. In particular, there exists  $\theta_0\in L^1(\Omega)$ such that $u_n\leq \theta_0$, that is, 
$$\|u_n^{\alpha+\gamma}\|_{q}\leq \left(|\Omega|+\|\theta_0\|_1\right)^{\frac{1}{q}}. $$

\nd If $(\alpha+\gamma)q>1$ we  distinguish between the two cases:  $\alpha> 1$ and $\alpha \leq 1$.
\vskip.1cm

\nd  In the case $\alpha> 1$, we find  by using that $((u_n+1/n)^{(\ell+\alpha-1)/\ell})$ is bounded in $W^{1,\ell}(\Omega)$ and
$W^{1,\ell}(\Omega)\hookrightarrow L^{\ell^*}(\Omega)$ that
$$\|u_n\|_{\ell^*+(\alpha-1)\frac{\ell^*}{\ell}}^{1+\frac{(\alpha-1)}{\ell}}=\left(\int_\Omega u_n^{\ell^*+(\alpha-1)\frac{\ell^*}{\ell}}dx\right)^{\frac{1}{\ell^*}}=\|u_n^{(\ell+\alpha-1)/\ell})\|_{\ell^*}\leq C,$$
that is, by using our assumption $q \leq q(\alpha+\gamma)$, it follows from its definition (see (\ref{probq}) for this) that $(\alpha+\gamma)q \leq \ell^*+(\alpha-1){\ell^*}/{\ell}$. So,
\begin{equation}\label{un-q-alpha-gamma}
\|u_n\|_{(\alpha+\gamma)q}\leq C,
\end{equation}
because 
 $L^{\ell^*+(\alpha-1){\ell^*}/{\ell}}(\Omega)\hookrightarrow L^{(\alpha+\gamma)q}(\Omega)$.
 \vskip.1cm
 
 \nd In the case $\alpha\leq 1$, again we have that $u_n$ is bounded in $\w$. So, it follows from $\w\hookrightarrow L^{(\gamma+\alpha)q}(\Omega)$, see (\ref{probq}) again, that 
$$\|u_n^{\alpha+\gamma}\|_{q}=\|u_n\|_{(\alpha+\gamma)q}^{\alpha+\gamma}\leq \kappa\|u_n\|^{\alpha+\gamma}\leq C,$$ for some $\kappa, C>0$. 
\vskip.1cm

\nd Thus, in both cases, it follows from (\ref{moser}) and the estimates just above  that  there exists a constant $c_0>0$ such that
\begin{eqnarray}\label{eq_1}
\frac{\beta_k}{q'}\int_\Omega\phi(|\nabla u_n|)|\nabla u_n|^2u_n^{\frac{\beta_k}{q'}+\alpha-1}dx&\leq&  (\|a\|_q+\|b\|_\infty c_0)\|u_n\|_{\beta_k}^{\frac{\beta_k}{q'}}.
\end{eqnarray}

\nd On the other hand, it follows by Lemma \ref{lema_naru} that
\begin{eqnarray}\label{eq_2}  
	\frac{\beta_k}{q'}\int_\Omega\phi(|\nabla u_n|)|\nabla u_n|^2u_n^{\frac{\beta_k}{q'}+\alpha-1}dx & \geq &\frac{\ell\Phi(1)}{q'}\beta_k\int_{|\nabla u_n|\geq 1}|\nabla u_n|^\ell u_n^{\frac{\beta_k}{q'}+\alpha-1} 
\end{eqnarray}
\nd and so it follows from  \eqref{eq_1} and  \eqref{eq_2}, that
\begin{eqnarray}\label{eq_3}
	\frac{\ell\Phi(1)}{q'}\beta_k\int_\Omega|\nabla u_n|^\ell u_n^{\frac{\beta_k}{q'}+\alpha-1}dx & \leq & \frac{\ell\Phi(1)}{q'}\beta_k\int_{|\nabla u_n|<1}|\nabla u_n|^\ell u_n^{\frac{\beta_k}{q'}+\alpha-1} dx \nonumber\\
	& + & (\|a\|_q+\|b\|_\infty c_0)\|u_n\|_{\beta_k}^{\frac{\beta_k}{q'}}\nonumber\\
	& \leq & \frac{\ell\Phi(1)}{q'}\beta_k\int_\Omega u_n^{\frac{\beta_k}{q'}+\alpha-1} dx \nonumber\\
	& + & (\|a\|_q+\|b\|_\infty c_0)\|u_n\|_{\beta_k}^{\frac{\beta_k}{q'}}.
\end{eqnarray}

\nd Our next objective is to show that 
\begin{eqnarray}\label{eq_5}
\int_\Omega|\nabla u_n|^\ell u_n^{\frac{\beta_k+(\alpha-1)q'}{q'}}dx & \leq &  B\|u_n\|_{\beta_k}^{\frac{\beta_k}{q'}},
\end{eqnarray}
\nd for some constant $B > 0$. To do this, we are going to consider two cases again: $\alpha \leq 1$ and $\alpha > 1$.
\vskip.2cm

\nd  If $\alpha \leq 1$ notice that  $L^{{\beta_k}}(\Omega)\hookrightarrow L^{\frac{\beta_k}{q'}+\alpha-1}(\Omega)$. Hence  
\begin{eqnarray}\label{eq-u_1}
	\int_\Omega u_n^{\frac{\beta_k}{q'}+\alpha-1} dx=\|u_n\|_{\beta_k/q'+\alpha-1}^{\beta_k/q'+\alpha-1}\leq |\Omega|^{1-\frac{1}{q'}+\frac{1-\alpha}{\beta_k}}\|u_n\|_{\beta_k}^{\beta_k/q'}\|u_n\|_{\beta_k}^{\alpha-1}.
\end{eqnarray}

\nd On the other hand, since $u_1\leq u_n$, we have 
\begin{equation}\label{eq-u_1-2}
	\|u_1\|_{\beta_k}\leq \|u_n\|_{\beta_k},
\end{equation} 
\nd and by the embedding
$L^{{\beta_k}}(\Omega)\hookrightarrow L^1(\Omega)$
\nd we get
\begin{equation}\label{eq-u_1-3}
	\|u_1\|_1\leq |\Omega|^{1-\frac{1}{\beta_k}}\|u_1\|_{\beta_k}.
\end{equation}
\nd Combining \eqref{eq-u_1-2} and \eqref{eq-u_1-3} we have
\begin{equation}\label{eq-u_1-4}
	\|u_n\|_{\beta_k}^{\alpha-1}\leq |\Omega|^{(1-\alpha)(1-\frac{1}{\beta_k})}\|u_1\|_1^{\alpha-1}.
\end{equation}

\nd So by \eqref{eq-u_1} and \eqref{eq-u_1-4} we infer that
\begin{eqnarray}\label{eq-u_1-5}
	\int_\Omega u_n^{\frac{\beta_k}{q'}+\alpha-1}dx\leq |\Omega|^{2-\alpha-\frac{1}{q'}}\|u_1\|_{1}^{\alpha-1}\|u_n\|_{\beta_k}^{\beta_k/q'}.
\end{eqnarray}

\nd Now by applying  \eqref{eq-u_1-5} in \eqref{eq_3}, we get

\begin{eqnarray}\label{eq_4}
	\frac{\ell\Phi(1)}{q'}\beta_k\int_\Omega|\nabla u_n|^\ell u_n^{\frac{\beta_k}{q'}+\alpha-1}dx & \leq & \frac{\ell\Phi(1)}{q'}|\Omega|^{2-\alpha-\frac{1}{q'}}\|u_1\|_{1}^{\alpha-1}\beta_k\|u_n\|_{\beta_k}^{\frac{\beta_k}{q'}} \nonumber\\
	& + & (\|a\|_q+\|b\|_\infty c_0)\|u_n\|_{\beta_k}^{\frac{\beta_k}{q'}}.
\end{eqnarray}

\nd Let $\alpha > 1$. By H\"older Inequality,  $(\alpha-1)q<(\alpha+\gamma)q$ and \eqref{un-q-alpha-gamma}, we have
\begin{eqnarray}\label{just-alpha-maior}
	\int_\Omega u_n^{\frac{\beta_k}{q'}+\alpha-1}dx&\leq& \|u_n\|_{\beta_k}^{\frac{\beta_k}{q'}}\left(\int_\Omega u_n^{(\alpha-1)q}dx\right)^{\frac 1 q}\nonumber\\
	&\leq &\|u_n\|_{\beta_k}^{\frac{\beta_k}{q'}}\left(|\Omega|+\int_{[u_n\geq 1]} u_n^{(\alpha-1)q}dx\right)^{\frac 1 q}\nonumber\\
	&\leq &\|u_n\|_{\beta_k}^{\frac{\beta_k}{q'}}\left(|\Omega|+\| u_n\|^{(\alpha+\gamma)q}_{(\alpha+\gamma)q}\right)^{\frac 1 q}\nonumber\\
	&\leq &\left(|\Omega|+C\right)^{\frac 1 q} \|u_n\|_{\beta_k}^{\frac{\beta_k}{q'}}.
\end{eqnarray} 
Now by applying \eqref{just-alpha-maior} in \eqref{eq_3}, we get
\begin{eqnarray}\label{eq_alph-maior}
\frac{\ell\Phi(1)}{q'}\beta_k\int_\Omega|\nabla u_n|^\ell u_n^{\frac{\beta_k}{q'}+\alpha-1}dx & \leq & \frac{\ell\Phi(1)}{q'}\beta_k\left(|\Omega|+C\right)^{\frac 1 q}\|u_n\|_{\beta_k}^{\frac{\beta_k}{q'}} \nonumber\\
& + & (\|a\|_q+\|b\|_\infty c_0)\|u_n\|_{\beta_k}^{\frac{\beta_k}{q'}}.
\end{eqnarray}

\nd So, it follows from \eqref{eq_4} (the case $\alpha \leq 1$) and \eqref{eq_alph-maior} (the case $\alpha > 1$) that  the inequality \eqref{eq_5} is true for $B>0$ defined by
$$
B:=\left\{
\begin{array}{l}
	\frac{q'}{\ell\Phi(1)}\left(\frac{\ell\Phi(1)}{q'}|\Omega|^{2-\alpha-\frac{1}{q'}}\|u_1\|_{1}^{\alpha-1}+\|a\|_q+\|b\|_\infty c_0\right), 0<\alpha\leq 1,\\
	\frac{q'}{\ell\Phi(1)}\left(\frac{\ell\Phi(1)}{q'}\left(|\Omega|+C\right)^{\frac 1 q}+\|a\|_q+\|b\|_\infty c_0\right), \alpha>1.
\end{array}
\right.
$$
This shows the inequality \eqref{eq_5}.
\nd Now since
\begin{eqnarray}\label{eq_6}
	\left(\frac{\ell q'}{\beta_k+\beta_1}\right)^\ell\int_\Omega\left|\nabla \left(u_n^{\frac{\beta_k+\beta_1}{\ell q'}}\right)\right|^\ell dx =\int_\Omega|\nabla u_n|^\ell u_n^{\frac{\beta_k+q'(\alpha-1)}{q'}}dx,\nonumber
\end{eqnarray}
it follows from 
 $\eqref{eq_5}$ and $W_0^{1,\ell}(\Omega)\hookrightarrow L^{\ell^*}(\Omega)$, that
\begin{eqnarray}\label{eq_8}
	\|u_n\|_{\beta_{k+1}}^{\frac{\beta_k^*}{q'}}=\left\|u^{\frac{\beta_k^*}{\ell q'}}\right\|_{\ell^*}^\ell  \leq \mu^\ell B\left(\frac{\beta_k^*}{\ell q'}\right)^\ell\|u_n\|_{\beta_k}^{\frac{\beta_k}{q'}},
\end{eqnarray}
for some $\mu>0$.

\nd Set $F_{k+1}:=\beta_{k+1}\ln (\|u_n\|_{\beta_{k+1}})$. So, it follows from the last inequality, that
\begin{eqnarray}
	F_{k+1} &\leq & \frac{\beta_{k+1}q'}{\beta_k^*}\left(\ell\ln \mu +\ell\ln\left(\frac{\beta_k^*}{\ell q'}\right)+\ln B+\frac{\beta_k}{q'}\ln(\|u_n\|_{\beta_k})\right)\nonumber\\
	&\leq& \ell^*\ln\left(\mu B\beta_k^*\right)+ \frac{\ell^*}{q'\ell}F_k\nonumber\\
	&=&\lambda_k+\delta F_k,
\end{eqnarray}
where $\lambda_k:=\ell^*\ln\left(\mu B\beta_k^*\right).$ 

Now, by using $(\ref{beta1})$ and $(\ref{beta2})$, we can infer that
$$\lambda_k=b+\ell^*\ln\left(2\delta^{k-1}+\delta^{k-2}+...+1\right),$$
where $b:=\ell^*\ln(\mu B\beta_1)$, that is, 
$$F_k\leq \delta^{k-1}F_1+\lambda_{k-1}+\delta\lambda_{k-2}+...+\delta^{k-2}\lambda_1.$$

So,
\begin{eqnarray}\label{eq_9}
 \frac{F_k}{\beta_k} & \leq & \frac{\delta^{k-1}F_1+\lambda_{k-1}+\delta\lambda_{k-2}+...+\delta^{k-2}\lambda_1}{\frac{2\delta^k-\delta^{k-1}-\delta}{\delta-1}\beta_1}\nonumber\\
 &=&\frac{F_1+\frac{\lambda_{k-1}}{\delta^{k-1}}+\frac{\lambda_{k-2}}{\delta^{k-2}}+...+\frac{\lambda_1}{\delta}}{\frac{2\delta-1-1/\delta^{k-1}}{\delta-1}\beta_1}.
\end{eqnarray}
Since 
\begin{eqnarray}\label{eq_10}
	\frac{\lambda_n}{\delta^n} &=& \frac{b}{\delta^n}+\frac{\ell^*}{\delta^n}\ln\left(\frac{2\delta^n-\delta^{n-1}-1}{\delta-1}\right)\nonumber\\
	&\leq& \frac{b}{\delta^n}+\frac{\ell^*}{\delta^n}\ln\left(\frac{2\delta^n}{\delta-1}\right),\nonumber
\end{eqnarray}
it follows from \eqref{eq_9}, that 
\begin{eqnarray}
	\frac{F_k}{\beta_k} & \leq & \frac{F_1+b\left(\frac{1}{\delta^{k-1}}+...\frac{1}{\delta}\right)+\ell^*\left(\frac{1}{\delta^{k-1}}\ln\left(\frac{2\delta^{k-1}}{\delta-1}\right)+...+
		\frac{1}{\delta}\ln\left(\frac{2\delta}{\delta-1}\right)\right)}{\frac{2\delta-1-1/\delta^{k-1}}{\delta-1}\beta_1}\nonumber\\
	&\leq & \frac{F_1+\frac{b}{\delta-1}+\ell^*\left(\frac{1}{\delta^{k-1}}\ln\left(\frac{2\delta^{k-1}}{\delta-1}\right)+...+
		\frac{1}{\delta}\ln\left(\frac{2\delta}{\delta-1}\right)\right)}{\frac{2\delta-1-1/\delta^{k-1}}{\delta-1}\beta_1}\nonumber\\
	& \leq & \frac{F_1+\frac{b}{\delta-1}+\ell^*\left[\ln\frac{2}{\delta-1}\left(\frac{1}{\delta^{k-1}}+...\frac{1}{\delta}\right)+\ln\delta\left(\frac{k-1}{\delta^{k-1}}+...\frac{1}{\delta}\right)\right]}{\frac{2\delta-1-1/\delta^{k-1}}{\delta-1}\beta_1}\nonumber\\
	& \leq & \frac{F_1+\frac{b}{\delta-1}+\ell^*\left[\frac{1}{\delta-1}\ln\frac{2}{\delta-1}+\ln\delta\sum_{n=1}^\infty \frac{n}{\delta^n}\right]}{\frac{2\delta-1-1/\delta^{k-1}}{\delta-1}\beta_1}\rightarrow d_0.
\end{eqnarray}

Now, going back to the definition of $F_k$, we obtain
$$\vert u_n(x)\vert \leq \|u_n\|_\infty=\limsup_{k\rightarrow\infty}\|u_n\|_{\beta_k}\leq \limsup_{k\rightarrow\infty} e^{\frac{F_k}{\beta_k}}\leq e^{d_0}~\mbox{for all}~ x\in \Omega,$$
and
$$\vert u(x)\vert = \lim_{n\to \infty}\vert u_n(x)\vert \leq e^{d_0}~\mbox{a.e}~ x\in \Omega,$$
because $u_n(x) \to u(x) $ a.e. in $\Omega$, that is, $u\in L^{\infty}(\Omega)$. This ends our proof. \hfill\cqd
\vskip.2cm

\nd{\bf Proof of Corollary 1.1:}  First, let $u \in \w$ be a solution of $(\ref{prob})$. Take $\varphi\in \w$, and $\varphi_n \in C^{\infty}_0(\Omega)$ such that $\varphi_n \to \varphi$ in $\w$. So, by taking $\sqrt[\theta]{\epsilon^{\theta}+ \vert \varphi_n - \varphi_k \vert^{\theta}} - \epsilon$, for some $\theta \in \mathbb{N}$, as a test function, we obtain
\begin{eqnarray}
0&\leq &\displaystyle \int_{\Omega}\frac{a(x)}{u^{\alpha}}\big[ \sqrt[\theta]{\epsilon^{\theta}+ \vert \varphi_n - \varphi_k \vert^{\theta}} - \epsilon\big]dx\\
&\leq&\displaystyle\int_{\Omega}\phi(\vert \nabla u \vert)\vert \nabla u \vert\frac{\vert \varphi_n - \varphi_k\vert^{\theta-1}}{\big[ {\epsilon^{\theta}+ \vert \varphi_n - \varphi_k \vert^{\theta}} \big]^{(\theta-1)/\theta}}\vert \nabla\varphi_n - \nabla\varphi_k \vert dx
\nonumber\\
&\leq &\Vert \phi(|\nabla u|)\vert \nabla u\vert \Vert_{L_{\tilde{\Phi}}(\Omega) }\Vert \nabla \varphi_n-\varphi_k \Vert_{L_{\Phi}(\Omega)} ,\nonumber
\end{eqnarray}
for every $\epsilon>0$ given.  Making $\epsilon \to 0$, we find that  
$$\Big(\frac{a \varphi_n}{u^{\alpha}}\Big)~\mbox{is a Cauchy sequence in}~ L^1(\Omega),$$
so that $({a \varphi_n})/{u^{\alpha}}\to \upsilon \in L^1(\Omega)$. Since, $\varphi_n (x) \to \varphi(x)$ a.e. in $\Omega$, we have that $\upsilon = ({a \varphi})/{u^{\alpha}}$. By hypothesis, $u\in \w$ satisfies $(see~(\ref{final1})
 )$
$$
\int_\Omega\phi(|\nabla u|)\nabla u\nabla \varphi_n dx=\int_\Omega\frac{a(x)}{u^\alpha}\varphi_n dx,
$$
\nd and,  passing to the limit, it follows that 
\begin{equation}
\label{fraca}
\int_\Omega\phi(|\nabla u|)\nabla u\nabla \varphi dx=\int_\Omega\frac{a(x)}{u^\alpha}\varphi dx~ \mbox{for all}~\varphi \in \w.
\end{equation}

\nd To complete the proof of the uniqueness, let $v\in \w$ be another solution of $(\ref{prob})$. Now assuming that  $u, v \in \w,~u \ne v$ satisfy $(\ref{fraca})$, setting  $\varphi=u-v$  and using  the fact that $\Delta_{\Phi}$ is stricly monotone, we obtain 
\begin{eqnarray}
0&\leq &\int_\Omega(\phi(|\nabla u|)\nabla u -\phi(|\nabla v|)\nabla v)(   \nabla u - \nabla v) dx\\
&=&\int_\Omega a(x)\Big(\frac{1}{u^\alpha}-\frac{1}{v^\alpha} \Big)(u-v)dx<0,\nonumber
\end{eqnarray}
\nd impossible. Now, we proceed to the regularity. First $(i)$. In this case, we have $a_n=a$ for $n$ large enough. So, as a consequence of the Comparison Principle, like at the end of the proof in Lemma $\ref{est}$, that $u_{n+1} \geq u_n$. Besides this, if we assume that
$$\Omega_0:=\Big\{x \in \Omega~/~ u_{n+1}(x)+\frac{1}{n+1}>u_{n}(x)+\frac{1}{n}\Big\}\subset\subset \Omega,$$
is not empty, then we would obtain 
$-\Delta_{\Phi}(u_{n+1}+1/(n+1)) \leq -\Delta_{\Phi}(u_{n}+1/n) $ in $\Omega_0$, that is, $ u_{n+1}(x)+\frac{1}{n+1}\leq u_{n}(x)+\frac{1}{n}$ in $\Omega_0$. This is impossible. So, we have
$$0 \leq u_n - u_k \leq \frac{1}{k} -\frac{1}{n}~\mbox{in}~ \Omega.$$
Since $(u_n)\subset C^1(\overline{\Omega})$, we obtain that $u_n$ converges uniformilly to $u$, that is, $u\in C(\overline{\Omega})$.

Proof of $(ii)$. It just follows from the same arguments that we used to proof Theorem $\ref{Teor-prin1}$, by taking  $b=0$. 
This ends our proof. \hfill\cqd

\section{Appendix  - On  Orlicz-Sobolev spaces}

\nd  In this section we present for, the reader's convenience, several results/notation used in the paper. 
 The reader is  referred to  $\cite{A,Rao1}$ regarding basics on Orlicz-Sobolev spaces.  The usual norm on $L_{\Phi}(\Omega)$ is, (Luxemburg norm),
\[
\|u\|_\Phi=\inf\left\{\lambda>0~|~\int_\Omega \Phi\left(\frac{u(x)}{\lambda}\right) dx \leq 1\right\},
\]
\nd while the  Orlicz-Sobolev norm of $ W^{1, \Phi}(\Omega)$ is
\[
\displaystyle \|u\|_{1,\Phi}=\|u\|_\Phi+\sum_{i=1}^N\left\|\frac{\partial u}{\partial x_i}\right\|_\Phi.
\]
\nd We denote by  $\w$  the closure of $C_0^{\infty}(\Omega)$ with respect to the Orlicz-Sobolev norm of $W^{1,\Phi}(\Omega)$. It remind that
$$
\widetilde{\Phi}(t) = \displaystyle \max_{s \geq 0} \{ts - \Phi(s) \},~ t \geq 0.
$$
\nd It turns out that  $\Phi$ and $\widetilde{\Phi}$  are  N-functions  satisfying  the $\Delta_2$-condition, (cf. \cite[p 22]{Rao1}).
In addition,   $L_{\Phi}(\Omega)$  and $W^{1,\Phi}(\Omega)$  are  reflexive and  Banach spaces.

\begin{rmk}\label{rem51}
	It is well known that $(\phi_3)$ implies that the condition
	\begin{itemize}
		\item[$(\phi_3)^\prime$]~~~~~~~~~~~~~~~~~ $\displaystyle\ell\leq \frac{\phi(t)t^2}{\Phi(t)}\leq m,~t>0,$
	\end{itemize}
	is verified.
	Furthermore, under this condition,  $\Phi,\widetilde\Phi\in\Delta_2$.
\end{rmk}

\nd By the Poincar\'e Inequality (see e.g.  \cite{gossez-Czech}), that is, the inequality
\[
\int_\Omega\Phi(u)dx\leq \int_\Omega\Phi(2d_{\Omega}|\nabla u|)dx,
\]
\nd where $d_{\Omega}=\mbox{diam}(\Omega)$,  it follows that
\[
\|u\|_\Phi\leq 2d_{\Omega}\|\nabla u\|_\Phi~\mbox{for all}~ u\in \w.
\]
\nd As a consequence of this, we have that  $\|u\| :=\|\nabla u\|_\Phi$ defines a norm in $\w$ that is equivalent to $\|\cdot\|_{1,\Phi}$. Let $\Phi_*$ be the inverse of the function
$$
t\in(0,\infty)\mapsto\int_0^t\frac{\Phi^{-1}(s)}{s^{\frac{N+1}{N}}}ds
$$
\nd which can be extended to ${\r}$ by  $\Phi_*(t)=\Phi_*(-t)$ for  $t\leq 0.$

We say that an N-function $\Psi$ grows essentially more slowly (grows more slowly) than $\Upsilon$, denoted by $\Psi<<\Upsilon$ $(\Psi<\Upsilon)$, if 
$$
\lim_{t\rightarrow \infty}\frac{\Psi(\lambda t)}{\Phi_*(t)}=0~~\mbox{for each}~~\lambda >0
$$
$(\Psi(t)\leq \Upsilon(kt)~\mbox{for all}~ t\geq t_0$ for some $k,t_0>0$).

The imbeddings below (cf. \cite{A}) were  used in this paper. First, we have
$$
\displaystyle \w \stackrel{\tiny cpt}\hookrightarrow L_\Psi(\Omega)~~\mbox{if}~~\Phi<\Psi<<\Phi_*,
$$
and in particular, 
$$
\w \stackrel{\tiny{cpt}} \hookrightarrow L_\Phi(\Omega),
$$
because $\Phi<<\Phi_*$ (cf. \cite[Lemma 4.14]{Gz1}). Furthermore,
$$
W_0^{1,\Phi}(\Omega) \stackrel{\mbox{\tiny cont}}{\hookrightarrow} L_{\Phi_*}(\Omega).
$$
Besides this, It is worth mentioning that if $(\phi_1)-(\phi_2)$ and $(\phi_3)^\prime$ are satisfied (cf. \cite[Lemma D.2]{clement}), then 
$$L_\Phi(\Omega)\stackrel{\mbox{\tiny cont}}\hookrightarrow L^\ell(\Omega).$$

We used in this text the notation $L^{\Psi}_{loc}(\Omega)$ in the sense that $u \in L^{\Psi}_{loc}(\Omega)$ if and only if $u \in L_{\Psi}(\Omega)$ for all $U\subset\subset \Omega$.

\begin{lem}\label{lema_naru}
{\rm (  cf. \cite{Fuk_1})}	Assume that  $\phi$ satisfies  $(\phi_1)-(\phi_3)$ hold.
	Set
	$$
	\zeta_0(t)=\min\{t^\ell,t^m\}~~\mbox{and}~~ \zeta_1(t)=\max\{t^\ell,t^m\},~~ t\geq 0.
	$$
	\nd Then  $\Phi$ satisfies
	$$
	\zeta_0(t)\Phi(\rho)\leq\Phi(\rho t)\leq \zeta_1(t)\Phi(\rho),~~ \rho, t> 0,
	$$
	$$
	\zeta_0(\|u\|_{\Phi})\leq\int_\Omega\Phi(u)dx\leq \zeta_1(\|u\|_{\Phi}),~ u\in L_{\Phi}(\Omega).
	$$
\end{lem}
\begin{lem}\label{lema_naru_*}
{\rm (  cf. \cite{Fuk_1})}	Assume that  $\phi$ satisfies $(\phi_1)-(\phi_3)$ and $1<\ell,m<N$ hold.  Set
	$$
	\zeta_2(t)=\min\{t^{\widetilde\ell},t^{\widetilde m}\}~~\mbox{and}~~ \zeta_3(t)=\max\{t^{\widetilde\ell},t^{\widetilde m}\},~~  t\geq 0,
	$$
	\nd where $\widetilde m = {m}/({m-1})$ and $\widetilde\ell = {\ell }/({\ell-1})$.  Then
	$$
	\widetilde\ell\leq\frac{t^2\widetilde\Phi'(t)}{\widetilde\Phi(t)}\leq \widetilde m,~t>0,
	$$
	$$
	\zeta_2(t)\widetilde\Phi(\rho)\leq\widetilde\Phi(\rho t)\leq \zeta_3(t)\widetilde\Phi(\rho),~~ \rho, t> 0,
	$$
	$$
	\zeta_2(\|u\|_{\widetilde\Phi})\leq\int_\Omega\widetilde\Phi(u)dx\leq \zeta_3(\|u\|_{\widetilde\Phi}),~ u\in L_{\widetilde\Phi}(\Omega).
	$$
\end{lem}

\begin{lem}\label{lema_Phi}
Let   Let $\Phi$ be an $N$-function satisfying  $\Delta_2$.  Let  $(u_n)\subset L_\Phi(\Omega)$ be a sequence such that $u_n\rightarrow u$ in  $L_\Phi(\Omega)$. Then there is a subsequence $(u_{n_k})\subseteq(u_n)$ such that
	\begin{description}
		\item{\rm{(i)}} $u_{n_k}(x)\rightarrow u(x) $ a.e. $x\in\Omega$,
		\item{\rm{(ii)}} there is  $h\in L_\Phi(\Omega)$ such that 
		$|u_{n_k}|\leq h~\mbox{a.e. in}~ \Omega.$
	\end{description}
\end{lem}

\nd \proof:(Sketch) We have that $\int_\Omega\Phi(u_n-u)dx\rightarrow 0$. By  \cite{A} $L_\Phi(\Omega)\hookrightarrow L^1(\Omega)$. So there are a  subsequence, we keep the notation, and  $\widetilde{h}\in L^1(\Omega)$ such that
$u_n\rightarrow u~\mbox{a.e. }~\mbox{in}~\Omega$ and $\Phi(u_n-u)\leq\widetilde{h}~\mbox{a.e.}~\mbox{in}~\Omega$.
\vskip.1cm

\nd Since  $\Phi$ is convex, increasing  and satisfies $\Delta_2$, we have
$$
\begin{array}{lll}
\Phi(|u_n|) & \leq & \displaystyle C\Phi\left(\frac{|u_n-u|+|u|}{2}\right) \leq  \displaystyle\frac{C}{2}\left[\Phi(|u_n-u|)+\Phi(|u|)\right]\\
& \leq& \displaystyle\frac{C}{2}[\widetilde{h}+\Phi(|u|)],
\end{array}
$$
that is, 
$$\vert u_n \vert \leq \Phi^{-1}\left(\frac{C}{2}(\widetilde{h}+\Phi(|u|))\right) :=h\in L_\Phi(\Omega),$$
because
 $\widetilde{h}\in L^1(\Omega)$, $\Phi(|u|)\in L^1(\Omega)$, and  
$$
\begin{array}{lll}
\displaystyle\int_\Omega\Phi(h)dx & = & \displaystyle\int_\Omega\Phi\left(\Phi^{-1}\left(\frac{K}{2}(\widetilde{h}+\Phi(|u|))\right)\right)dx \\ \\
& = & \displaystyle \int_\Omega\left(\frac{K}{2}(\widetilde{h}+\Phi(|u|))\right)dx <  \infty,
\end{array}
$$
\nd showing that $h\in L_\Phi(\Omega)$. \hfill\cqd

{\tiny }

\begin{flushright}
\scriptsize{ Marcos L. M. Carvalho \\

 Jos\'e V.A. Gon\c{c}alves   }\\
\smallskip
  \scriptsize{Universidade Federal de Goi\'as\\
   Instituto de Matem\'atica e Estat\'istica\\
   74001-970 Goi\^ania, GO - Brasil}

 \scriptsize{ emails: marcosleandrocarvalho@yahoo.com.br\\

goncalves.jva@gmail.com}
\end{flushright}

\begin{flushright}
{\scriptsize Carlos Alberto Santos}\\
	\smallskip
	\scriptsize{Universidade de Bras\'ilia\\
		Departamento de Matem\'atica\\
		Bras\'ilia, DF 70910--900, Brazil}\\

\scriptsize{ email:  csantos@unb.br}\\

\end{flushright}

\end{document}